\documentclass[10pt]{article}%
\usepackage{amsfonts}
\usepackage[left=3.2cm,right=3.2cm,top=2cm,bottom=4cm]{geometry}%
\usepackage{amsmath}%
\usepackage{pstricks}
\usepackage{amssymb}%
\usepackage{graphicx}
\usepackage{multirow}
\usepackage{enumitem}
\usepackage{amsthm}
\usepackage{pdflscape}
\usepackage{listings}

\providecommand{\U}[1]{\protect\rule{.1in}{.1in}}
\newcommand{\esssup}{\operatorname*{ess\;sup}}

\newtheorem{theorem}{Theorem}

\newtheorem{lemma}[theorem]{Lemma}

\newtheorem{proposition}[theorem]{Proposition}
\newtheorem{remark}[theorem]{Remark}

\newcommand{\sgn}{\operatorname{sgn}}
\newcommand{\rec}{(\ref{eq:recA}),(\ref{eq:recB}),(\ref{eq:recinvA}) and (\ref{eq:recinvB}) }

\newcommand{\N}{\mathbb{N}}
\newcommand{\R}{\mathbb{R}}

\begin{document}
\lstset{language=Mathematica}

\title{The Lebesgue Constant for the Periodic Franklin System}
\author{Markus Passenbrunner\footnote{Supported by FWF P 20166-N18.}
\footnote{
This is part of my PhD thesis written at Department of Analysis, J. Kepler University Linz. I want to thank my advisor P. F. X M\"uller for many helpful discussions during the preparation of this thesis.}}
\maketitle
\renewcommand{\thefootnote}{}
\tableofcontents
\newpage

\begin{abstract}
We identify the torus with the unit interval $[0,1)$ and let $n,\nu\in\mathbb{N}$, $0\leq \nu\leq n-1$ and $N:=n+\nu$. Then we define the (partially equally spaced) knots
\[
t_{j}=\left\{
\begin{array}
[c]{ll}%
\frac{j}{2n}, & \text{for }j=0,\ldots,2\nu,\\
\frac{j-\nu}{n}, & \text{for }j=2\nu+1,\ldots,N-1.
\end{array}
\right. 
\]
Furthermore, given $n,\nu$ we let $V_{n,\nu}$ be the space of piecewise linear continuous functions on the torus with knots $\left\{t_j:0\leq j\leq N-1\right\}$. Finally, let $P_{n,\nu}$ be the orthogonal projection operator from $L^{2}([0,1))$ onto
$V_{n,\nu}.$ The main result is
\[
\lim_{n\rightarrow\infty,\nu=1}\left\|P_{n,\nu}:L^\infty\rightarrow L^\infty\right\|=\sup_{n\in\mathbb{N},0\leq \nu\leq n}\left\|P_{n,\nu}:L^\infty\rightarrow L^\infty\right\|=2+\frac{33-18\sqrt{3}}{13}.
\]
This shows in particular that the Lebesgue constant of the classical Franklin orthonormal system on the torus is $2+\frac{33-18\sqrt{3}}{13}$.
\end{abstract}

\section{Introduction}
Let $(N_k)_{k\geq 0}$ be an orthonormal basis in $L^2[0,1]$. The Fourier partial sums with respect to this basis are given by
\begin{equation}
P_N(f)=\sum_{k=0}^N \left<f,N_k\right>N_k. \label{eq:foursum}
\end{equation}
Clearly, every $P_N$ is a projection onto its (finite dimensional) range and its norm as an operator from $L^\infty[0,1]$ to $L^\infty[0,1]$ (or as an operator from $L^1[0,1]$ to $L^1[0,1]$) is given by
\[
L_N=\esssup_{s\in [0,1]}\int_0^1 |K_N(s,t)|dt,
\]
where $K_N$ is the \emph{Dirichlet kernel}
\[
K_N(s,t)=\sum_{k=0}^N N_k(s)N_k(t).
\]
The \emph{Lebesgue constant} of the basis $(N_k)_{k\geq 0}$ is now defined as
\[
L:=\sup_{N\geq 0} L_N.
\]
As a particular instance of an orthonormal basis in $L^2[0,1]$, we consider the \emph{general Franklin system} $(N_k)_{k\geq 0}$ on the torus $\mathbb{T}=\R/\mathbb{Z}$: That is we choose a sequence of points $\mathcal{T}=(t_k)_{k\geq 0}$ in $[0,1)$ (we identify this interval with the torus), which is dense in $[0,1)$ and with $t_0=0$. The space of piecewise linear and continuous functions on $\mathbb{T}$ with knots $\{t_0,\ldots, t_N\}$ is denoted by $V_N(\mathcal{T})$. Then we define $f_0\equiv 1$ on $\mathbb{T}$ and inductively, for $k\geq 1$ the $k$-th Franklin function corresponding to the sequence $\mathcal{T}$ is uniquely determined by the conditions
\[
f_k\in V_k(\mathcal{T}),\quad f_k \perp V_{k-1}(\mathcal{T}),\quad \left\|f_k\right\|_2=1,\quad f_k(t_k)>0.
\]
The Franklin functions $f_k$ are splines of degree $d=1$. We now make a few comments about the history of calculating or estimating the Lebesgue constant of splines of degree $d$.

For $d=0$ (piecewise constant functions), the projection is easily calculated and the Lebesgue constant is $1$.

For $d=1$ (piecewise linear functions), Z. Ciesielski (\cite{Ciesielski1963}) proved that for any partition $\pi$ of $[0,1]$, the $L^\infty$-norm onto piecewise linear functions with knots $\pi$ is $\leq 3$. He showed this for the non-periodic case, but exactly the same argument gives the upper bound $3$ in the periodic case. Moreover, P. Oswald (\cite{Oswald1977}) and K. Oskolkov (\cite{Oskolkov1979}) proved independently that in the non-periodic case, the constant $3$ is optimal if one considers arbitrary partitions $\pi$. Moreover, Ciesielski (\cite{Ciesielski1975a}) showed that in case of uniform partitions the exact upper bound is $2$. Some numerical experiments suggested that for the (classical, corresponding to dyadic knots) non-periodic Franklin system, the  exact upper bound is $2+(2-\sqrt{3})^2$ (\cite{CiesielskiNiedzwiecka}). Several years later, P. Bechler (\cite{Bechler2003}) proved that for the piecewise linear Str\"{o}mberg wavelet, the Lebesgue constant is indeed $2+(2-\sqrt{3})^2$. Then, Z. Ciesielski and A. Kamont (\cite{CiesielskiKamont2004}) showed that for the classical non-periodic Franklin system, the Lebesgue constant is $2+(2-\sqrt{3})^2$, verifying the conjecture in \cite{CiesielskiNiedzwiecka}.

For splines of higher degree ($d\geq 2$), a problem was the mere existence of a bound $C_d$ for the $L^\infty-$norms of orthogonal projections onto splines of degree $d$ with arbitrary knots, where $C_d$ depends only on $d$ and not on the partition. This was a long standing conjecture by C. de Boor solved by A. Yu Shadrin in \cite{Shadrin2001} (in the non-periodic case). Predating Shadrin's result, there were several results specializing in the degree (for instance \cite{deBoor1968} for $d=2$ in the non-periodic case) or specializing in the sequence of points (for instance \cite{Domsta1972} and \cite{Domsta1976} viewing the sequence of dyadic partitions both in the non-periodic and periodic case respectively for arbitrary degree $d$). 
In the periodic case, there is a further partial result showing the existence of a bound $C_2$ for the $L^\infty-$norm of orthogonal projections for $d=2$ not depending on the knots in \cite{Keryan2008}. The exact values of the Lebesgue constants in the cases $d\geq 2$ are not known.

In the present paper, we study and determine the Lebesgue constant for the periodic (classical) Franklin system (corresponding to $d=1$). Its value is $2+\frac{33-18\sqrt{3}}{13}$. The analysis presented in this article was constantly guided by extensive computer simulations (both numerically and symbolically) involving the Gram matrix and its inverse (see Section \ref{sec:orth}).

\paragraph{Acknowledgements} I am grateful to A. Kamont and the anonymous referee who made many valuable comments and suggestions to earlier versions of this article.

\section{Formulation of the Main Theorem}
Our main result concerns partially equally spaced knots on the torus $\mathbb{T}=\mathbb{R}/\mathbb{Z}$. We choose the special points
\begin{equation}
t_{j}=\left\{
\begin{array}
[c]{ll}%
\frac{j}{2n}, & \text{for }j=0,\ldots,2\nu\\
\frac{j-\nu}{n}, & \text{for }j=2\nu+1,\ldots,N-1
\end{array}
\right.  \label{eq:specpoints}
\end{equation}
for arbitrary $n,\nu\in\mathbb{N}$ with $0\leq \nu\leq n-1$ and $N:=n+\nu$.  We remark that for $\nu=0$ or $\nu=n$ we arrive at equally spaced knots. Let $V_{n,\nu}$ be the linear subspace generated by the piecewise linear, continuous functions with knots (\ref{eq:specpoints}) and $P_{n,\nu}$ be the orthogonal projection onto $V_{n,\nu}$. The B-spline basis for $V_{n,\nu}$ with a special choice of parameters $n,\nu$ is pictured in Figure \ref{fig:splines}.

The main theorem now reads as follows:
\begin{theorem}\label{th:main}
For all $n\in\mathbb{N}, 0\leq\nu\leq n$, we have the following bound for the norm of the projection operator $P_{n,\nu}$ onto $V_{n,\nu}$:
\[
\left\|P_{n,\nu}\right\|_\infty:=\left\| P_{n,\nu}:L^\infty(\mathbb{T})\rightarrow L^\infty(\mathbb{T})\right\|< 2+\frac{33-18\sqrt{3}}{13}=:\gamma.
\]
Furthermore, for $n\rightarrow\infty, \nu=1$ it holds that
\[
\lim_{n\rightarrow\infty}\left\|P_{n,1}\right\|_\infty=\gamma.
\]
\end{theorem}

\section{Preliminaries}
\subsection{Orthogonal Projections}\label{sec:orth}
Let $V$ be an $N$-dimensional subspace of $L^2[0,1]$ and $\{N_0,\ldots,N_{N-1}\}$ a basis of $V$. We first look at the changes in formula (\ref{eq:foursum}), if the basis functions are no longer orthogonal. In this case, the orthogonal projection $P$ onto $V$ is given by
\[
P f(s)=\sum_{j,k=0}^{N-1} a_{jk}\left<N_k,f\right>N_j(s),
\]
or equivalently as an integral operator with kernel $k(s,t)=\sum_{j,k=0}^{N-1} a_{jk}N_j(s)N_k(t)$
\[
P f(s)=\int_0^1 k(s,t)f(t)dt,
\]
where $(a_{jk})$ is the inverse of the Gram matrix $(b_{jk})$ with $b_{jk}=\left<N_j,N_k\right>$. The norm of $P$ as a mapping from $L^\infty [0,1]$ to $L^\infty [0,1]$ is
\begin{equation}
\left\|P\right\|_\infty=\esssup_{s\in [0,1]} \int_0^1 |k(s,t)|dt.\label{eq:projnorm}
\end{equation}
Since $P$ is self adjoint, the norm of $P$ as operator from $L^1[0,1]$ to $L^1[0,1]$ is the same.

We now consider periodic B-splines of degree one on $\mathbb{T}=\mathbb{R}/\mathbb{Z}$. For this let $0=t_0<t_1<\cdots <t_{N-1}<1$ with an arbitrary natural number $N\geq 2$. Further set $t_{-1}:=t_{N-1}-1$, $t_N:=1$ and $\delta_j:=t_{j+1}-t_j$ for $-1\leq j\leq N-1$. Then we let $N_j$ for $0\leq j\leq N-1$ be the unique continuous function on $\mathbb{T}$, which is linear on every interval $(t_{k-1},t_k)$ and has values $N_j(t_k)=\delta_{j,k}$ for $0\leq k\leq N-1$. Formally we define the functions $N_j:\mathbb{T}\rightarrow [0,1]$ for $0\leq j\leq N-1$ as 
\begin{equation}
N_j([t]):=\begin{cases}
(s-t_{j-1})/\delta_{j-1}, &  \text{if } [t]=[s]\text{ for } t_{j-1}<s\leq t_j, \\
(t_{j+1}-s)/\delta_j, & \text{if } [t]=[s]\text{ for } t_j<s\leq t_{j+1}, \\
0, & \text{otherwise,}
\end{cases} \label{eq:bsplines}
\end{equation}
where we denote by $[\cdot]$ the canonical surjection taking each $t\in\mathbb{R}$ onto its equivalence class in $\mathbb{T}$. From now on we identify the unit interval $[0,1)$ with $\mathbb{T}$ and furthermore, by a slight abuse of notation we consider $N_j$ to be defined on $[0,1)$. 

\begin{figure}[h]
\begin{center}
\begin{pspicture}(-2.5,-0.3)(11.25,3.5)
\psline(-2.5,0)(11.25,0)
\psline(0,0)(0,3)
\psline(10,0)(10,3)
\psline(0,-0.15)(0,0.15)
\psline(5,-0.15)(5,0.15)
\psline(2.5,-0.15)(2.5,0.15)
\psline(7.5,-0.15)(7.5,0.15)
\psline(1.25,-0.15)(1.25,0.15)
\psline(10,-0.15)(10,0.15)
\psset{linewidth=1.5pt}
\psline(-2.5,0)(0,3)(1.25,0)
\psline(0,0)(1.25,3)(2.5,0)
\psline(1.25,0)(2.5,3)(5,0)
\psline(2.5,0)(5,3)(7.5,0)
\psline(5,0)(7.5,3)(10,0)
\psline(7.5,0)(10,3)(11.25,0)
\uput[d](0,-0.15){$t_0$}
\uput[d](1.25,-0.15){$t_1$}
\uput[d](2.5,-0.15){$t_2$}
\uput[d](5,-0.15){$t_3$}
\uput[d](7.5,-0.15){$t_4$}
\uput[d](10,-0.15){$t_0$}
\uput[d](0.625,0){$\delta_0$}
\uput[d](1.875,0){$\delta_1$}
\uput[d](3.75,0){$\delta_2$}
\uput[d](6.25,0){$\delta_3$}
\uput[d](8.75,0){$\delta_4=\delta_{-1}$}
\uput[u](0,3){$N_0$}
\uput[u](1.25,3){$N_1$}
\uput[u](2.5,3){$N_2$}
\uput[u](5,3){$N_3$}
\uput[u](7.5,3){$N_4$}
\uput[u](10,3){$N_0$}
\end{pspicture}
\end{center}
\caption{Situation for $N=5,\nu=1,n=N-\nu=4$.}
\label{fig:splines}
\end{figure}
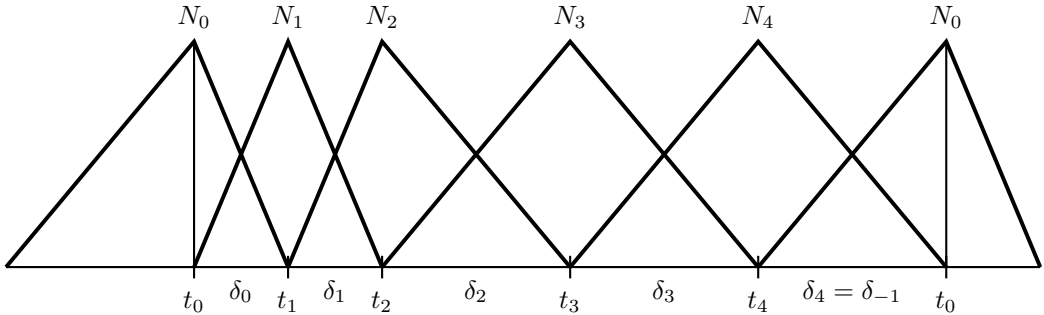

Figure \ref{fig:splines} shows periodic B-splines of degree one defined in (\ref{eq:bsplines}) for the points in (\ref{eq:specpoints}) with a special choice of parameters $n,\nu$.

Let (as above) $V$ be the (finite dimensional) subspace generated by $\{N_0,\ldots, N_{N-1}\}$ and $P$ be the orthogonal projection from $L^2[0,1)$ onto $V$. Then formula (\ref{eq:projnorm}) for the norm of $P$ simplifies to
\[
\left\|P\right\|_\infty=\max_{j=0,\ldots,N-1}\int_0^1 |k(t_j,t)|dt,
\]
where the kernel $k$ is given by $k(s,t)=\sum_{j,k=0}^{N-1} a_{j,k}N_j(s)N_k(t)$. Now recall that $(a_{j,k})$ is the inverse of the Gram matrix $(b_{j,k})=\left<N_j,N_k\right>$. If we let $\kappa(j):=\int_0^1 |k(t_j,t)|dt$, it can be shown by an elementary calculation that
\begin{equation}
\label{eq:kappainvgram}
\kappa(j)=\sum_{k=0}^{N-1}\frac{\delta_k}{2}
\begin{cases}
|a_{j,k}|+|a_{j,k+1}|, & \text{ if }\sgn a_{j,k}=\sgn a_{j,k+1}, \\
\frac{a_{j,k}^2+a_{j,k+1}^2}{|a_{j,k}|+|a_{j,k+1}|}, & \text{ otherwise},
\end{cases}
\end{equation}
where every subindex is understood to be an index modulo $N$. Observe that $\kappa(j)$ depends on $N$ too.
With the rational function $\phi(t):=\frac{1+t^2}{(1+t)^2}$, equation $(\ref{eq:kappainvgram})$ can be rewritten to
\begin{equation}\label{eq:formulaprojnorm}
\kappa(j)=\sum_{k=0}^{N-1}\frac{\delta_k}{2}(|a_{j,k}|+|a_{j,k+1}|)\cdot
\begin{cases}
1, & \text{ if }\sgn a_{j,k}=\sgn a_{j,k+1}, \\
\phi(|a_{j,k+1}|/|a_{j,k}|), & \text{ otherwise}.
\end{cases}
\end{equation}
We now collect a few simple facts about the function $\phi$:
\begin{lemma}\label{lem:phi}
Let $\phi:(0,\infty)\rightarrow [1/2,1)$ be defined by
\[
t\mapsto\phi(t)=\frac{1+t^2}{(1+t)^2}.
\]
Then
\[
\phi(t)=\phi(t^{-1}),\quad\phi'(t)=\displaystyle\frac{2(t-1)}{(1+t)^3},\quad\phi''(t)=\displaystyle\frac{4(2-t)}{(1+t)^4}
\]
for all $t> 0$. So in particular $\phi$ is decreasing for $t<1$ and increasing for $t>1$ and $\phi'$ is increasing for $t<2$ and decreasing for $t>2$. Furthermore,
\[
\phi(\lambda)=\frac{2}{3},\quad\phi(4)=\frac{17}{25},\quad\phi(6)=\frac{37}{49},\quad \phi'(\lambda)=\frac{\lambda^{-1}}{3\sqrt{3}},
\]
where $\lambda=2+\sqrt{3}$.
\end{lemma}

By (\ref{eq:formulaprojnorm}), exact formulae for the entries of the inverse $(a_{jk})$ of the Gram matrix are absolutely necessary in determining the exact value of the Lebesgue constant. We will provide this information in Proposition \ref{prop:formelgj0} for the periodic case.
In the non-periodic dyadic case, such exact formulae for the inverse of the Gram matrix were given in \cite{Ciesielski1966} and they were used in the calculation of the corresponding Lebesgue constant in \cite{CiesielskiKamont2004}. For the general Franklin system, there are important estimates both for the non-periodic case and for the periodic case (see \cite{KashinSaakyan1989} and \cite{Keryan2005} respectively). To calculate the exact value of the Lebesgue constant, we supplemented these already known estimates with exact formulae. 


\subsection{Solutions of $f_{k-1}-4f_{k}+f_{k+1}=0$ and their Properties}
In this section we define and examine a few properties of the solutions of the recurrence $f_{k-1}-4f_{k}+f_{k+1}=0$, which we will use extensively in the sequel.
For an arbitrary real number $x$, let $A_{x}:=\cosh(\alpha x)$ and $\sqrt{3}B_{x}:=\sinh(\alpha x)$ with
$\alpha>0$ defined by $\cosh\alpha=2.$ For $k\in\N_0$, $A_{k}$ and $B_{k}$ can also be defined
by the recurrence relations%
\begin{eqnarray}
A_{k+1}&=& 2A_{k}+3B_{k}\quad\text{with }A_{0}=1,\label{eq:recA}\\
B_{k+1}&=& A_{k}+2B_{k}\quad\text{with }B_{0}=0.\label{eq:recB}
\end{eqnarray}
This follows from the basic identities
\begin{eqnarray}
\cosh(x+y)  & =&\cosh x\cosh y+\sinh x\sinh y,\label{eq:hyp1}\\
\sinh(x+y)  & =&\sinh x\cosh y+\cosh x\sinh y.\label{eq:hyp2}%
\end{eqnarray}
We note that it is easy to see (or a special case of Lemma \ref{lem:asym}) that the inequalities 
\begin{eqnarray}
A_{k+1}&\leq &4A_k\quad\text{for }k\in\N_0, \label{eq:A4}\\
B_{k+1}&\leq &4B_k\quad\text{for }k\in\N \label{eq:B4}
\end{eqnarray}
hold. Observe also that%
\begin{align}
A_{k}  & =2A_{k+1}-3B_{k+1}\label{eq:recinvA}\\
B_{k}  & =2B_{k+1}-A_{k+1} \label{eq:recinvB}
\end{align}
for $k\in\mathbb{N}_0$. We also have the formulae%
\begin{equation}
A_{x}=\frac{1}{2}(\lambda^{x}+\lambda^{-x}),\quad B_{x}=\frac{1}{2\sqrt{3}}(\lambda
^{x}-\lambda^{-x}),\quad x\in\R \label{eq:defAB}
\end{equation}
with%
\[
\lambda=2+\sqrt{3},\quad \lambda^{-1}=2-\sqrt{3}.
\]
We remark that $\alpha=\log\lambda$.
For reference, we list the first few values of both $A_{n}$ and $B_{n}:$%
\begin{align*}
(A_{0},\ldots,A_{4}) =(1,2,7,26,97),\qquad
(B_{0},\ldots,B_{4}) =(0,1,4,15,56).
\end{align*}

The crucial fact about $A_k$ and $B_k$ is that they are independent solutions of the linear recursion $f_{k-1}-4f_{k}+f_{k+1}=0$, since $\lambda$ and $\lambda^{-1}$ are the two solutions of its characteristic equation $t^2-4t+1=0$ and $A_k$ and $B_k$ have the representation (\ref{eq:defAB}). The recursion $f_{k-1}-4f_{k}+f_{k+1}=0$ in turn takes into account the special form of the Gram matrix for the points (\ref{eq:specpoints}) (see (\ref{eq:grameq}) and (\ref{eq:gramineq})). This is important, since we need exact formulae for the inverse of the Gram matrix and these consist then of terms depending on $A_k$ and $B_k$.

\begin{lemma}
\label{lem:sumid}For $K\in\mathbb{N}_0$ we have the following formulae%
\begin{eqnarray*}
\sum_{k=0}^{K}B_{k}+B_{k+1}=A_{K+1}-1, \quad 2\sum_{k=0}^K A_k=3B_{K+1}-A_{K+1}+1, \\
\sum_{k=0}^{K}A_{k}+A_{k+1}=3B_{K+1}, \quad 2\sum_{k=0}^K B_k=A_{K+1}-B_{K+1}-1.
\end{eqnarray*}
\end{lemma}

\begin{proof}
The proof uses induction and the recurrences (\ref{eq:recA}),(\ref{eq:recB}),(\ref{eq:recinvA}) and (\ref{eq:recinvB}) for $A_n$ and $B_n$.
\end{proof}

\begin{lemma}
\label{lem:asym}Let $k\in\mathbb{N}_0.$ Then we have 
\begin{eqnarray}
-1  \leq -\lambda^{-k}&=&\lambda B_{k}-B_{k+1}\leq 0, \\
0\leq\lambda A_{k}-A_{k+1}&=&\sqrt{3}\lambda^{-k}\leq \sqrt{3},\\
-1  \leq\lambda^{-k}&=&\sqrt{3}B_{k}-A_{k}\leq0.
\end{eqnarray}

\end{lemma}

\begin{proof}
This follows from (\ref{eq:defAB}).
\end{proof}

\begin{lemma}
\label{lem:trigid}For all $n\in\mathbb{N}$ and $0\leq k\leq n$ the following
equalities hold
\begin{eqnarray*}
B_kA_{n-k}+A_kB_{n-k}=B_n, \quad B_nA_{n-k}-B_{n-k}A_n=B_k,   \\
A_kA_{n-k}+3B_{n-k}B_k=A_n, \quad A_nA_{n-k}-3B_nB_{n-k}=A_k.
\end{eqnarray*}
\end{lemma}
\begin{proof}
This follows directly from (\ref{eq:hyp1}) and (\ref{eq:hyp2}).
\end{proof}

\section{Proof of the Main Theorem}
We begin with a short overview of the main steps of the proof. In Section \ref{sec:eqspaced} we treat the special case of equally spaced knots, since this is the simplest one and we get an even better Lebesgue constant than the one stated in Theorem \ref{th:main} here. This serves as some kind of preliminary result, where all important proof-steps of more general cases are included:
\begin{enumerate}
\item Compute the inverse of the Gram matrix.
\item Estimate $L^\infty$-norms of the projection operators using step 1. For this, it is important to distinguish the cases where the number of points in the knot sequence is even or odd. This difference in the analysis comes from the fact that the inverse of the Gram matrix has a different structure depending on this distinction. 
\item Determine the asymptotics of these projection operator norms.
\end{enumerate}

In Section \ref{sec:invgram} we calculate the inverse of the Gram matrix for non-equally spaced knots. \\
Section \ref{sect:nu1j1} concentrates on estimating $\left\|P_{n,\nu}\right\|_\infty$ for the parameter choice $\nu=1$, where (as we will see) we get the largest values for the projection operator norms. We furthermore determine the asymptotics in this case which gives us the asserted value $2+\frac{33-18\sqrt{3}}{13}$ of the Lebesgue constant. \\
In Section 4.4 we estimate the remaining cases for other parameter choices of $\nu$ by employing easy-to-use, but sufficiently sharp 
estimates on quotients of consecutive entries of the inverse of the Gram matrix. 

\subsection{Equally Spaced Knots}\label{sec:eqspaced}
As a preliminary case we view the points (\ref{eq:specpoints}) for $\nu=0$ and $N=n$ and show that the $L^\infty$-norm $\|P_{n,0}\|_\infty$ obeys the estimate $\|P_{n,0}\|_\infty<2$ and that $\lim_{n\rightarrow\infty}\|P_{n,0}\|_\infty=2$. For this case of equally spaced knots, the Gram matrix $(b_{jk})_{0\leq j,k\leq N-1}$ is
\begin{equation}
(b_{jk})=\frac{1}{6n}\left(\begin{tabular}{p{0.47cm}p{0.47cm}p{0.47cm}p{0.47cm}c}
  4 & 1 & & & 1 \\
  1 & 4 & 1 &  &  \\
  & $\ddots$ & $\ddots$ & $\ddots$ &  \\  
  & & 1 & 4 & 1 \\  
  1 & & & 1 & 4
 \end{tabular}\right),
 \label{eq:grameq}
\end{equation}
where the empty entries are zero. Since every row in $(b_{jk})$ is equal up to shifts, the same must be true for the inverse $(a_{jk})$. For the first row of $(a_{jk})$, make the ansatz $a_{0,k}=(-1)^k(c_1 A_k+c_2 B_k)$ with constants $c_1,c_2$ that are to be determined. Thus it holds that
\[
a_{0,k}+4a_{0,k+1}+a_{0,k+2}=0\quad\text{for }k\geq 0
\] 
Insert this ansatz into the boundary conditions
\[
4a_{0,0}+a_{0,1}+a_{0,N-1}=1,\quad a_{0,N-2}+4a_{0,N-1}+a_{0,0}=0
\]
to determine $c_1,c_2$ and simplify to get
\[
a_{0,k}=\frac{6n (-1)^k}{D(N)}g_k
\]
with 
\begin{equation}
g_k=B_{N-k}+(-1)^NB_k\quad\text{and}\quad D(N)=2((-1)^{N-1}+A_N). 
\end{equation}
Since every row in $(a_{jk})$ is equal up to shifts, formula (\ref{eq:formulaprojnorm}) does not depend on $j$ in this case. So while we consider equally spaced knots, we write $\kappa$ to denote the value of $\kappa(j)$ for arbitrary $0\leq j\leq N-1$.
We consider separately the cases $N$ even and $N$ odd. The difference in the analysis of these two cases comes from the fact that $g_k$ is always positive for $N$ even, whereas for $N$ odd the sign of $g_k$ changes once.
\paragraph{$N$ even}

If we let $N$ even, we obtain from (\ref{eq:formulaprojnorm})
\begin{eqnarray*}
\kappa&=&3D(N)^{-1}\sum_{k=0}^{N-1} (g_k+g_{k+1})\phi\left(\frac{g_{k+1}}{g_k}\right).
\end{eqnarray*}
Using the definition of $g_k$ and Lemma \ref{lem:asym} we see that $\lambda^{-1}< \frac{g_{k+1}}{g_k}< \lambda$, so by Lemma \ref{lem:phi}, $\phi\left(\frac{g_{k+1}}{g_k}\right)<\phi(\lambda)$ and thus we obtain 
\[
\kappa< 6\phi(\lambda)D(N)^{-1}\sum_{k=0}^{N-1} B_k+B_{k+1}.
\]
Lemma \ref{lem:sumid} and the fact that $\phi(\lambda)=\frac{2}{3}$ then give us
\[
\kappa< 4\frac{A_N-1}{2(A_N-1)}=2.
\]

\paragraph{$N$ odd}

For $N$ odd, we see that (\ref{eq:formulaprojnorm}) becomes
\begin{eqnarray*}
\kappa&=&6D(N)^{-1} \Bigg[ B_{(N+1)/2}-B_{(N-1)/2}+\\ 
&& \sum_{j=0}^{(N-3)/2}(B_{N-j}+B_{N-j-1}-B_j-B_{j+1})\phi\left(\frac{B_{N-j}-B_j}{B_{N-j-1}-B_{j+1}}\right)\Bigg].
\end{eqnarray*}
The mean value theorem implies
\[
\phi(q_j)\leq \phi(\lambda)+(q_j-\lambda)\phi'(\lambda),\quad\text{where}\quad q_j:=\frac{B_{N-j}-B_j}{B_{N-j-1}-B_{j+1}},
\]
since $\phi'(t)$ is decreasing for $t\geq \lambda\geq 2$ and $q_j\geq \lambda$ by Lemma \ref{lem:asym}. For $q_j-\lambda$, we have again due to Lemma \ref{lem:asym} and $0\leq j\leq (N-3)/2$
\begin{eqnarray*}
q_j-\lambda&=&\frac{B_{N-j}-\lambda B_{N-j-1}+\lambda B_{j+1}-B_j}{B_{N-j-1}-B_{j+1}}\leq \frac{1+\lambda B_{j+1}}{B_{N-j-1}-B_{j+1}}\\
&\leq& \frac{1+\lambda B_{j+1}}{B_{N-j-1}(1-\lambda^{-N+2j+2})}\leq 2\frac{1+\lambda B_{j+1}}{B_{N-j-1}}.
\end{eqnarray*}
If we use these facts and the estimates $B_{(N-1)/2}\geq \lambda^{-1} B_{(N+1)/2}-\lambda^{-1}$ (Lemma \ref{lem:asym}) and $-B_j\leq 0$, we obtain for $\kappa$
\begin{eqnarray}
\kappa&\leq& 6D(N)^{-1}\Bigg[(1-\lambda^{-1})B_{(N+1)/2}+\lambda^{-1}+ \nonumber\\
&& +\phi(\lambda)\sum_{j=0}^{(N-3)/2} (B_{N-j}+B_{N-j-1}-B_{j+1})\nonumber \\
&& +2\phi'(\lambda)\sum_{j=0}^{(N-3)/2} (B_{N-j}+B_{N-j-1}-B_{j+1})\frac{1+\lambda B_{j+1}}{B_{N-j-1}}\Bigg]\label{eq:kappagleich}
\end{eqnarray}

We split the analysis of this expression into a few subcases and thereby introduce the notation $p=\frac{N+1}{2}$ to shorten indices.
\begin{enumerate}[labelindent=\parindent,leftmargin=*,label=\textsc{Sum }\Roman*.,widest=III,align=left]
\item $\sum_{j=0}^{p-2} B_{N-j}+B_{N-j-1}-B_{j+1}$ \\
We apply Lemma \ref{lem:sumid} and get that
\begin{eqnarray*}
\sum_{j=0}^{p-2} B_{N-j}+B_{N-j-1}-B_{j+1}&=&\frac{1}{2}(2A_N-3A_{p}+B_{p}+1) \\
&\leq & \frac{1}{2}(2A_N-(3\sqrt{3}-1)B_{p} +1),
\end{eqnarray*}
by Lemma \ref{lem:asym}.
\item $II:=\sum_{j=0}^{p-2}(B_{N-j}+B_{N-j-1}-B_{j+1})\frac{1+\lambda B_{j+1}}{B_{N-j-1}}$ \\
Since by Lemma \ref{lem:asym}, $B_{N-j}=\lambda B_{N-j-1}+\lambda^{-N+j+1}$ and $\lambda^{-N+j+1}\leq \lambda^{-N+1}B_{j+1}$, we get that
\[
II\leq(1+\lambda)\sum_{j=0}^{p-2}(1+\lambda B_{j+1})-(1-\lambda^{-N+1})\sum_{j=0}^{p-2}\frac{B_{j+1}(1+\lambda B_{j+1})}{B_{N-j-1}}.
\]
But now, by estimating the second sum by its summand with index $p-2$
\[
\sum_{j=0}^{p-2}\frac{B_{j+1}(1+\lambda B_{j+1})}{B_{N-j-1}}\geq \frac{B_{p-1}(1+\lambda B_{p-1})}{B_{p}}\geq\lambda\frac{B_{p-1}^2}{B_{p}},
\]
and by Lemmas \ref{lem:sumid} and \ref{lem:asym}
\[
\sum_{j=0}^{p-2}1+\lambda B_{j+1}=\frac{N-1}{2}+\frac{\lambda}{2}(A_{p}-B_{p}-1)\leq \frac{N-1}{2}+\frac{\lambda (\sqrt{3}-1)}{2}B_{p}
\]
We thus obtain finally
\[
II\leq (1+\lambda)\left(\frac{N-1}{2}+\frac{\lambda (\sqrt{3}-1)}{2}B_{p}\right)-(1-\lambda^{-N+1})\lambda\frac{B_{p-1}^2}{B_{p}}.
\]
\end{enumerate}
These estimates and (\ref{eq:kappagleich}) yield, noting $D(N)\geq 2A_N$ and $\phi(\lambda)=2/3$,
\[
\kappa\leq 2+\frac{3}{A_N}\left[\theta B_p+\lambda^{-1}+\frac{1}{3}+2\phi'(\lambda)\left((1+\lambda)\frac{N-1}{2}-(1-\lambda^{-N+1})\lambda\frac{B_{p-1}^2}{B_p}\right)\right],
\]
where
\[
\theta=(1-\lambda^{-1})-(\sqrt{3}-\frac{1}{3})+(1+\lambda)\lambda\phi'(\lambda)(\sqrt{3}-1)=0.
\]
Since $\frac{B_{p-1}^2}{B_p}$ dominates $(N-1)/2$ for large $N$, we finally get that for $N$ sufficiently large ($N\geq 8$)
\[
\kappa<2.
\]
In fact, if we look at Table \ref{tab:numvals} on page \pageref{tab:numvals}, we see that for all $N\geq 2$ we have this inequality.
An analogous argument as in Section \ref{sec:asym} finally yields that $\lim_{N\rightarrow\infty}\kappa=2$, and this completes what we wanted to show in this section.
\subsection{The Inverse of the Gram Matrix for Non-Equally Spaced Knots}\label{sec:invgram}
We now view the points (\ref{eq:specpoints}) in case $1\leq \nu\leq n-1$ (i.e. the case where the knots are not equally spaced anymore). The first step is to calculate the inverse of the Gram matrix in this setting, which we do in this section. As above and in the following we understand every index concerning the Gram matrix $(b_{jk})$ or its inverse $(a_{jk})$ as an index modulo $N$. The Gram matrix $(b_{jk})=(\left\langle N_{j},N_{k}\right\rangle)_{0\leq j,k\leq N-1}$ admits the following representation
\begin{equation}
(b_{jk})=\frac{1}{12n}\left(\begin{tabular}{p{0.47cm}p{0.47cm}p{0.47cm}p{0.47cm}p{0.47cm}p{0.47cm}p{0.47cm}p{0.47cm}c}
  6 & 1 & & & & & & & 2 \\
  1 & 4 & 1 & & & & & &  \\
  & $\ddots$ & $\ddots$ & $\ddots$ & & & & \\
  & & 1 & 4 & 1 & & & &\\ 
  & & & 1 & 6 & 2 & & & \\
  & & & & 2 & 8 & 2 & & \\
  & & & & & $\ddots$ & $\ddots$ & $\ddots$ & \\
  & & & & & & 2 & 8 & 2 \\
  2 & & & &  & & & 2 & 8 \\  
 \end{tabular}\right),
 \label{eq:gramineq}
\end{equation}
where the row with the pattern $1, 6, 2$ has the index $2\nu.$
This leads to the following equations concerning the inverse $(a_{jk})$ of $(b_{jk})$:
\begin{align}
6a_{0,k}+a_{1,k}+2a_{N-1,k}  & =12n\delta_{0,k},\label{eqn:1}\\
a_{j-1,k}+4a_{j,k}+a_{j+1,k}  & =12n\delta_{j,k}\quad\text{for }%
j=1,\ldots,2\nu-1,\label{eqn:2}\\
a_{2\nu-1,k}+6a_{2\nu,k}+2a_{2\nu+1,k}  & =12n\delta_{2\nu,k},\label{eqn:3}\\
a_{j-1,k}+4a_{j,k}+a_{j+1,k}  & =6n\delta_{j,k}\quad\text{for }j=2\nu
+1,\ldots,N-1,\label{eqn:4}
\end{align}
where $\delta_{j,k}$ is the Kronecker delta and $0\leq k\leq N-1$.
Let 
\begin{equation}
D(N,\nu):=2A_{N}+\frac{3}{2}B_{2\nu}B_{N-2\nu}-2(-1)^N.\label{eq:D}
\end{equation}
Then we define
\[
g(N,\nu,j,k):=\frac{D(N,\nu)a_{j,k}(-1)^{k+j}}{6n}
\]
Observe that $a_{j,k}$ depends on $N$ and $\nu$ too. But in the current context, the indices $N$,$\nu$ and also $j$ are fixed, so we write $g_k$ instead of $g(N,\nu,j,k)$. Inserting the definition of $g_k$ into (\ref{eq:formulaprojnorm}), we obtain
\begin{equation}\label{eq:kappa}
\kappa(j)=D(N,\nu)^{-1}\left[\frac{3}{2}\sum_{k=0}^{2\nu-1}(|g_k|+|g_{k+1}|)\cdot \xi_{j,k}+3\sum_{k=2\nu}^{N-1}(|g_k|+|g_{k+1}|)\cdot \xi_{j,k}\right]
\end{equation}
with
\[
\xi_{j,k}=\begin{cases}
1, & \text{if }\sgn a_{j,k}=\sgn a_{j,k+1},\\
\phi(|g_{k+1}|/|g_k|),& \text{else}.
\end{cases}
\]


In order to determine $(a_{j,k})$, we identify the values of $g_k$. This is the content of the following
\begin{proposition}\label{prop:formelgj0}
If $0\leq j\leq 2\nu-1$ we have that $g_k$ equals
\[
\begin{array}{ll}
2(-1)^N B_{j-k}+B_{N-j+k}+B_{2\nu-j}A_{N-2\nu+k}+B_k(A_{N-j}+3B_{2\nu-j}B_{N-2\nu}), & \text{if } k\leq j, \\
2(-1)^N B_{k-j}+B_{N-k+j}+B_{2\nu-k}A_{N-2\nu+j}+B_j(A_{N-k}+3B_{2\nu-k}B_{N-2\nu}), & \text{if }j\leq k\leq 2\nu,\\
(-1)^N(B_{k-j}+A_{k-2\nu}B_{2\nu-j})+B_{N-k+j}+B_jA_{N-k}, & \text{if }2\nu\leq k\leq N-1.	
\end{array}
\]
If $2\nu\leq j\leq N-1$, we have that $g_k$ equals
\[
\begin{array}{ll}
(-1)^N(B_{j-k}+A_{j-2\nu}B_{2\nu-k})+B_{N-j+k}+A_{N-j}B_k, & \text{if } k\leq 2\nu\leq j, \\
(-1)^N B_{j-k}+A_{k-2\nu}B_{N-j+2\nu}+A_{N-j}B_k+\frac{3}{2}B_{k-2\nu}B_{2\nu}B_{N-j},& \text{if } 2\nu \leq k\leq j, \\
(-1)^N B_{k-j}+A_{N-k}B_j+A_{j-2\nu}B_{N-k+2\nu}+\frac{3}{2}B_{2\nu}B_{N-k}B_{j-2\nu}, & \text{if } j\leq k\leq N-1.
\end{array}
\]
\end{proposition}

\begin{proof}
If we insert these formulae for $g_k$ into equations (\ref{eqn:1}) and (\ref{eqn:2}) for $0\leq j\leq 2\nu-1$ and into equations (\ref{eqn:3}) and (\ref{eqn:4}) for $2\nu\leq j\leq N-1$, we see the assertion of the proposition after a few case distinctions and uses of the fact that $A_n$ and $B_n$ are solutions of the recurrence $f_{k-1}-4f_k+f_{k+1}=0$. Observe that for evaluating (\ref{eqn:1}),(\ref{eqn:2}),(\ref{eqn:3}),(\ref{eqn:4}) the recursions (\ref{eq:recA}),(\ref{eq:recB}),(\ref{eq:recinvA}),(\ref{eq:recinvB}) for $A_k$ and $B_k$ and the identities from Lemma \ref{lem:trigid} are useful.
\end{proof}
\begin{remark}\label{eq:remgpos}
From the formulae in Proposition \ref{prop:formelgj0} 
we obtain that for $N$ even, $g_k\geq 0$ for all $0\leq k\leq N-1$ and for $N$ odd it holds that $g_k\geq 0$ for $|k-j|\leq \frac{N-1}{2}$ and $g_k\leq 0$ for $|k-j|\geq \frac{N+1}{2}$.
\end{remark}
\subsection{The Main Case $\nu=j=1$}\label{sect:nu1j1}
The first special case to analyze is the parameter choice $\nu=j=1$. As we will see in the sequel, this is the main case in the sense that for $N\rightarrow\infty$ and $\nu=j=1$, $\kappa:=\kappa(1)$ converges to the Lebesgue constant $2+\frac{33-18\sqrt{3}}{13}$.
In this section, we set $K=N-1$ for notational convenience. We then get as a special instance of Proposition \ref{prop:formelgj0}
\[
g_k=g(N,1,1,k)=\left\{
\begin{array}
[c]{ll}%
2\left[(-1)^N+A_{K}-B_{K}\right], & \text{if }k=0,\\
8B_{K}, & \text{if }k=1,\\
2\left[A_{N-k}+B_{N-k}+(-1)^N(A_{k-2}+B_{k-2})\right], & \text{if }2\leq k\leq N-1.
\end{array}
\right.
\]
Note that $g_2=g_0.$ Additionally it holds that%
\[
D(N,1)=18B_{K}-2A_{K}-2(-1)^N.
\]
Furthermore the use of the recurrences \rec for $A_k$ and $B_k$ yields
\begin{eqnarray}
|g_1|+|g_2|&=&2(-1)^N+6B_{K}+2A_{K}, \label{eq:sumg1}\\
|g_k|+|g_{k+1}|&=&4|A_{N-k}+(-1)^NA_{k-1}|\quad\text{for }k\geq 2,k\neq (N+1)/2, \\
|g_{(N+1)/2}|+|g_{(N+3)/2}|&=&8A_{K/2}\quad\text{for }N \text{ even},\\
|g_{(N+1)/2}|+|g_{(N+3)/2}|&=&8B_{K/2}\quad\text{for }N \text{ odd}\label{eq:sumg2}.
\end{eqnarray}
We recall that all indices have to be taken modulo $N$. The quotient of subsequent values of $g_k$ has the following special form
\begin{lemma}\label{lem:eqquot}
For $2\leq k\leq N-1$ it holds that
\begin{eqnarray}\label{eq:quotg1}
\frac{|g_{k+1}|}{|g_k|}&=&\frac{A_{|N/2-k|}}{A_{|N/2-k+1|}}\quad\text{for }N \text{ even,}\\
\label{eq:quotg2}
\frac{|g_{k+1}|}{|g_k|}&=&\frac{B_{|N/2-k|}}{B_{|N/2-k+1|}}\quad\text{for }N \text{ odd.}
\end{eqnarray}
\end{lemma}
\begin{proof}
Let $k\leq N/2$. Then we have by (\ref{eq:hyp1}), (\ref{eq:hyp2}) and the definitions of $A_n$ and $B_n$
\begin{eqnarray*}
A_{N-k-1}&=&A_{N/2-k}A_{N/2-1}+3B_{N/2-k}B_{N/2-1}, \\
B_{N-k-1}&=&A_{N/2-k}B_{N/2-1}+B_{N/2-k}A_{N/2-1}, \\
A_{k-1}&=&A_{N/2-1}A_{N/2-k}-3B_{N/2-k}B_{N/2-1}, \\
B_{k-1}&=&B_{N/2-1}A_{N/2-k}-A_{N/2-1}B_{N/2-k}.
\end{eqnarray*}
For $N$ even, summing these four equations yields $g_{k+1}/2$ on the left hand side and $A_{N/2-k}$ times a term independent of $k$ on the right hand side. On the other hand, for $N$ odd, summing the first two equations and subtracting the second two gives us $|g_{k+1}|/2$ on the left hand side and $B_{N/2-k}$ times a term independent of $k$ on the right hand side.
An analogous argument for $k\geq N/2$ completes the proof of the lemma.
\end{proof}

\subsubsection{Estimates for $N$ even}
For $N$ even, we get from (\ref{eq:kappa}) and the fact $g_0=g_2$ that
\[
\kappa:=\kappa(1)=3D(N,1)^{-1}\sum_{k=1}^{K}(g_k+g_{k+1})\phi\left(\frac{g_{k+1}}{g_k}\right).
\]
Inserting (\ref{eq:sumg1})-(\ref{eq:sumg2}) into this expression for $\kappa$ and recalling $K=N-1$ yield that $\kappa$ equals
\begin{eqnarray}
 3D(N,1)^{-1}\left[(2+6B_{K}+2A_{K})\phi\left(\frac{1+A_{K}-B_{K}}{4B_{K}}\right)+8\sum_{k=2}^{K}A_{k-1}\phi\left(\frac{A_{|N/2-k|}}{A_{|N/2-k+1|}}\right)\right]\label{eq:kappagerade}.
\end{eqnarray}
Now observe that Lemma \ref{lem:eqquot} and Lemma \ref{lem:asym} imply $\lambda^{-1}<\frac{g_{k+1}}{g_k}=\frac{A_{|N/2-k|}}{A_{|N/2-k+1|}}< \lambda$ for $k\geq 2$, so using Lemma \ref{lem:phi}, the previous expression for $\kappa$ is strictly less than
\[
6D(N,1)^{-1}\left[(1+3B_{K}+A_{K})\phi\left(\frac{1+A_{K}-B_{K}}{4B_{K}}\right)+4\phi(\lambda)\sum_{k=1}^{K-1}A_k\right].
\]
If we use Lemma \ref {lem:sumid} to evaluate the sum and remark that 
$A_{K}=\sqrt{3}B_{K}+\lambda^{-K}$ by Lemma \ref{lem:asym}, we obtain by setting $\frac{1+A_{K}-B_{K}}{4B_{K}}=\eta+h$ with
\[
\eta=\frac{\sqrt{3}-1}{4}\quad \text{and}\quad h=h(N)=\frac{1+\lambda^{-K}}{4B_{K}}
\]
the subsequent estimate for $\kappa$:
\begin{equation}\label{eq:genau1}
\kappa\leq 6D(N,1)^{-1}\left[(1+3B_{K}+A_{K})\phi(\eta+h)+2\phi(\lambda)(3B_{K}-A_{K}-1)\right].
\end{equation}
Since $\phi'(t)$ is increasing for $t\leq 2$ (Lemma \ref{lem:phi}) and $h\leq 1/2$ for $N\geq 2$, we get by applying the mean value theorem to $\phi$:
\begin{equation}\label{eq:genau2}
\phi(\eta+h)\leq \phi(\eta)+\phi'(\eta+\frac{1}{2})h.
\end{equation}
Thus, using (\ref{eq:genau2}) in (\ref{eq:genau1}) we see that in order to prove $\kappa< \gamma$, it suffices to show that
\begin{equation}
\label{eq:genau3}
6D(N,1)^{-1}\left[(1+3B_{K}+A_{K})(\phi(\eta)+\phi'(\eta+1/2)h)+2\phi(\lambda)(3B_{K}-A_{K}-1)\right]<\gamma.
\end{equation}
If we multiply this inequality by $D(N,1)$, collect the factors for $B_K$ and $A_K$ and observe that 
\[
\theta:=6\phi(\eta)+2\gamma-12\phi(\lambda)=\frac{1}{\sqrt{3}}(18\gamma-18\phi(\eta)-36\phi(\lambda)),
\]
we see that (\ref{eq:genau3}) is equivalent to
\begin{equation}
\label{eq:genau4}
\theta(\sqrt{3}B_K-A_K-1)+6h(N)(1+3B_K+A_K)|\phi'(\eta+1/2)|> 0.
\end{equation}
Now we use again $A_{K}=\sqrt{3}B_{K}+\lambda^{-K}$ and insert the definition of $h(N)$ to express the left hand side of (\ref{eq:genau4}) as
\[
(1+\lambda^{-K})\left[\frac{3}{2B_K}(1+(\sqrt{3}+3)B_K+\lambda^{-K})|\phi'(\eta+1/2)|-\theta\right].
\]
Clearly, this is greater than
\[
(1+\lambda^{-K})\left[\frac{3(\sqrt{3}+3)}{2}|\phi'(\eta+1/2)|-\theta\right]
\]
and this is easily seen to be greater than zero. Thus we have shown for $N$ even and $\nu=j=1$ that $\kappa< \gamma.$

\subsubsection{Estimates for $N$ odd}
For $N$ odd, (\ref{eq:kappa}) and Remark \ref{eq:remgpos} yield for $\kappa$ the formula
\[
\kappa=3D(N,1)^{-1}\left[\sum_{\stackrel{k=1}{k\neq (N+1)/2}}^{K}(|g_k|+|g_{k+1}|)\phi\left(\frac{|g_k|}{|g_{k+1}|}\right)+|g_{(N+1)/2}|+|g_{(N+3)/2}|\right].
\]
We now use Lemma \ref{lem:eqquot} and the identities (\ref{eq:sumg1})-(\ref{eq:sumg2}) and recall the setting $K=N-1$ to obtain after a little calculation that
\begin{eqnarray}
\kappa&=&6D(N,1)^{-1}\Bigg[(3B_{K}+A_{K}-1)\phi\left(\frac{A_K-B_K-1}{4B_K}\right) \nonumber\\
&& +4\sum_{k=2}^{K/2}(A_{N-k}-A_{k-1})\phi\left(\frac{B_{N/2-k}}{B_{N/2-k+1}}\right)+4B_{K/2}\Bigg]\label{eq:kappaungerade}.
\end{eqnarray}
We first estimate two summands of $\kappa$ separately
\setlist{noitemsep}
\setlength{\parindent}{0cm}
\begin{enumerate}[labelindent=\parindent,leftmargin=*,label=\textsc{Term }\Roman*.,widest=III,align=left]
\item $(3B_{K}+A_{K}-1)\phi\left(\frac{A_K-B_K-1}{4B_K}\right)$. \\
We have $3B_K+A_K-1\leq (3+\sqrt{3})B_K$ by Lemma \ref{lem:asym} and $\frac{A_K-B_K-1}{4B_K}=\eta-h$ with 
\[
\eta=\frac{\sqrt{3}-1}{4}\quad \text{and}\quad h=\frac{1-\lambda^{-K}}{4B_K},
\]
so the mean value theorem implies
\begin{eqnarray*}
(3B_{K}+A_{K}-1)\phi\left(\frac{A_K-B_K-1}{4B_K}\right)&\leq& (3+\sqrt{3})B_K \phi(\eta-h) \\
&\leq& (3+\sqrt{3})B_K (\phi(\eta)-\phi'(0)h)\\
&=&(3+\sqrt{3})B_K (\phi(\eta)+2h),
\end{eqnarray*}
since $\phi'$ is increasing for $t\leq 2$ and $\phi'(0)=-2$.

\item $II:=\sum_{k=2}^{K/2}(A_{N-k}-A_{k-1})\phi\left(\frac{B_{N/2-k+1}}{B_{N/2-k}}\right)$. \\
Since $B_{L+1}=\lambda B_L+\lambda^{-L}$, we get with the mean value theorem and the fact that $\phi'$ is decreasing for $t\geq 2$ 
\[
\phi\left(\frac{B_{N/2-k+1}}{B_{N/2-k}}\right)\leq \phi(\lambda)+\phi'(\lambda)\frac{\lambda^{k-N/2}}{B_{N/2-k}}.
\]
Now, if we use the identity $2\sum_{k=0}^L A_k=3B_{L+1}-A_{L+1}+1$ from Lemma \ref{lem:sumid} and simplify using the recurrences for $A_k$ and $B_k$, we obtain 
\begin{eqnarray*}
\sum_{k=2}^{K/2}A_{N-k}-A_{k-1}&=&\frac{1}{2}(3B_{K}-A_K-6B_{K/2}+1)\\
&\leq&\frac{1}{2}((3-\sqrt{3})B_{K}-6B_{K/2}+1),
\end{eqnarray*}
by Lemma \ref{lem:asym}. Next, we get
\begin{eqnarray*}
S:=\sum_{k=2}^{K/2} A_{N-k}\frac{\lambda^{k-N/2}}{B_{N/2-k}}&=&\sqrt{3}\sum_{k=2}^{K/2}\frac{\lambda^{k-N/2}(\lambda^{N-k}+\lambda^{k-N})}{\lambda^{N/2-k}-\lambda^{k-N/2}}\\
&=&\sqrt{3}\sum_{k=2}^{K/2}\frac{\lambda^{N-k}+\lambda^{k-N}}{\lambda^{N-2k}-1},
\end{eqnarray*}
by (\ref{eq:defAB}). Since $1\leq \lambda^{N-2k}/2$, we estimate
\begin{eqnarray*}
S&\leq& 2\sqrt{3}\sum_{k=2}^{K/2}\frac{\lambda^{N-k}+\lambda^{k-N}}{\lambda^{N-2k}}= 2\sqrt{3}\sum_{k=2}^{K/2}\lambda^{k}+\lambda^{3k-2N}\\
&=&2\sqrt{3}\left[\frac{\lambda^{K/2+1}-\lambda^2}{\lambda-1}+\lambda^{-2N}\frac{\lambda^{3(K/2+1)}-\lambda^6}{\lambda^3-1}\right]\\
&\leq& 2\sqrt{3}\left[\frac{\lambda^{K/2+1}}{\lambda-1}+\lambda^{-2N}\frac{\lambda^{3(K/2+1)}}{\lambda-1}\right] \\
&=&4\sqrt{3}\frac{A_{K/2}}{1-\lambda^{-1}}\leq 4\sqrt{3}\frac{\sqrt{3}B_{K/2}+1}{1-\lambda^{-1}}.
\end{eqnarray*}
If we summarize all estimates, we get for the whole sum
\begin{eqnarray*}
II\leq\frac{\phi(\lambda)}{2}((3-\sqrt{3})B_K-6B_{K/2}+1)+4\sqrt{3}\phi'(\lambda)\frac{\sqrt{3}B_{K/2}+1}{1-\lambda^{-1}}.
\end{eqnarray*}
\end{enumerate}
Let us now return to (\ref{eq:kappaungerade}). If we use the estimate $h\leq \frac{1}{4B_K}$, we obtain by combining the estimates for the Terms I and II that
\[
\kappa\leq 6D(N,\nu)^{-1}(\sigma B_K-\tau B_{K/2}+\vartheta),
\]
with $\sigma=(3+\sqrt{3})\phi(\eta)+2\phi(\lambda)(3-\sqrt{3})$, $\tau=12\phi(\lambda)-4-\frac{48\phi'(\lambda)}{1-\lambda^{-1}}>0$ and $\vartheta=\frac{3+\sqrt{3}}{2}+2\phi(\lambda)+\frac{16\sqrt{3}\phi'(\lambda)}{1-\lambda^{-1}}$. Now recall that $D(N,1)=18B_K-2A_K+2\geq (18-2\sqrt{3})B_K$ by Lemma \ref{lem:asym}, so in order to prove $\kappa<\gamma$, it suffices to show
\[
\frac{\gamma (18-2\sqrt{3})B_K}{6}> \sigma B_K-\tau B_{K/2}+\vartheta.
\]
Since $\sigma=\frac{\gamma}{6}(18-2\sqrt{3})$, this is equivalent to
\[
\tau B_{K/2}-\vartheta>0,
\]
which is true for $N\geq 7$. For $N<7$ we get the desired bound for $\kappa$ from Table \ref{tab:numvals} on page \pageref{tab:numvals}.
\subsubsection{Asymptotic Behaviour}\label{sec:asym}
In this section, we calculate the limit of $\kappa$ as $N\rightarrow \infty$ for $\nu=j=1$. In the following, the symbol $\sim$ will denote asymptotic equality for $N\rightarrow\infty$. If we remark $A_N\sim\sqrt{3}B_N$, $A_{N+1}\sim\lambda A_N$ (by Lemma \ref{lem:asym}) and recall the definition of $D(N,1)=18B_K-2A_K-2$ (where as above, $K=N-1$) we get for $N$ even from (\ref{eq:kappagerade})
\begin{eqnarray*}
\kappa&=& 6D(N,1)^{-1}\left[(1+3B_{K}+A_{K})\phi\left(\frac{1+A_{K}-B_{K}}{4B_{K}}\right)+4\sum_{k=2}^{K}A_{k-1}\phi\left(\frac{A_{|N/2-k|}}{A_{|N/2-k+1|}}\right)\right] \\
&\sim& \frac{6}{(18-2\sqrt{3})B_K}\left[(3+\sqrt{3})B_K\phi\left(\frac{\sqrt{3}-1}{4}\right)+4\phi(\lambda)\sum_{k=3N/4}^{K-1} A_{k}\right].
\end{eqnarray*}
Using the identity $2\sum_{k=0}^L A_k=3B_{L+1}-A_{L+1}+1$ from Lemma \ref{lem:sumid}, we get further
\begin{eqnarray*}
\kappa&\sim& \frac{6}{(18-2\sqrt{3})B_K}\left[(3+\sqrt{3})B_K\phi\left(\frac{\sqrt{3}-1}{4}\right)+2\phi(\lambda)(3-\sqrt{3})B_K\right] \\
&\sim& \frac{6}{18-2\sqrt{3}}\left[(3+\sqrt{3})\phi\left(\frac{\sqrt{3}-1}{4}\right)+2\phi(\lambda)(3-\sqrt{3})\right]\\
&=&\gamma=2+\frac{33-18\sqrt{3}}{13}.
\end{eqnarray*}
If on the other hand $N$ is odd, we obtain from (\ref{eq:kappaungerade})
\begin{eqnarray*}
\kappa
&\sim& \frac{6}{(18-2\sqrt{3})B_K}\left[(3+\sqrt{3})B_K\phi\left(\frac{\sqrt{3}-1}{4}\right)+4\sum_{k=2}^{K/4}A_{N-k}\phi\left(\frac{B_{N/2-k}}{B_{N/2-k+1}}\right)\right].
\end{eqnarray*}
Again, the identity $2\sum_{k=0}^L A_k=3B_{L+1}-A_{L+1}+1$ and $B_{N+1}\sim\lambda B_N$ imply in the same way as above
\[
\kappa\sim \gamma.
\]
Thus if we combine the estimates of this section (Section \ref{sect:nu1j1}) with the numerical results from Table \ref{tab:numvals} on page \pageref{tab:numvals} we have shown that for $\nu=j=1$, we have $\kappa<\gamma$ and $\lim_{N\rightarrow\infty}\kappa=\gamma$. We will see in the next section, that this is the critical case, since we will show that for all other values of $\nu$ and $j$ we have $\kappa<\gamma$.
\subsection{Estimating $\kappa(j)$}\label{sec:estkappa}
In this section we derive bounds for $\kappa(j)$ for all remaining values of $\nu,j$, which will allow us to deduce that for all $n,\nu\in\N$, $0\leq \nu\leq n$, we have $\|P_{n,\nu}\|_\infty<\gamma$.
In order to derive these estimates for $\kappa(j)$ we first need some for the quotients of subsequent values of $g$. This is the content of the following two lemmas.
\begin{lemma}\label{lem:bdg}
Let $N$ be even. Then it holds that 
\begin{eqnarray}
6^{-1}&\leq& \frac{g_{k+1}}{g_k}\leq 6\quad\text{for }k=0 \text{ or } k=2\nu-1, \label{eq:schrg1}\\
4^{-1}&\leq& \frac{g_{k+1}}{g_k}\leq 4\quad\text{for }k\neq 0\text{ and } k\neq 2\nu-1. \label{eq:schrg2}
\end{eqnarray}
For $j=0,k=0$, we have a better estimate
\[
4^{-1}\leq \frac{g_{k+1}}{g_k}\leq 4.
\]
\end{lemma}

We get analogous estimates for $N$ odd, but we have to add a further restriction to the domain of validity of the inequalities:
\begin{lemma}\label{lem:schrunger}
Let $N\geq 7$ be odd and $|k-j|\leq\frac{N-5}{2}$ or $|k-j|\geq \frac{N+5}{2}$. Then we have
\begin{eqnarray*}
6^{-1}&\leq& \frac{|g_{k+1}|}{|g_k|}\leq 6 \quad\text{for }k=0\text{ or } k=2\nu-1, \\
4^{-1}&\leq& \frac{|g_{k+1}|}{|g_k|}\leq 4 \quad\text{for }k\neq 0\text{ and } k\neq 2\nu-1.
\end{eqnarray*}
Additionally, for $j=0,k=0$ we have the better estimate
\[
4^{-1}\leq \frac{|g_{k+1}|}{|g_k|}\leq 4.
\]
\end{lemma}
For a proof of Lemma \ref{lem:bdg} or parts of a proof of Lemma \ref{lem:schrunger}, see Appendix A. 

We note that in the following, we only treat the case $N$ even. In fact, as we will show later (in Section \ref{sec:kappaodd}), the case $N$ odd will nonetheless follow from these estimates. Combining formula (\ref{eq:kappa}) with Remark \ref{eq:remgpos} yields for $N$ even
\begin{equation}
\kappa(j)=D(N,\nu)^{-1}\left[\frac{3}{2}\sum_{k=0}^{2\nu-1}(g_k+g_{k+1})\phi\left(\frac{g_{k+1}}{g_k}\right)+3\sum_{k=2\nu}^{N-1}(g_k+g_{k+1})\phi\left(\frac{g_{k+1}}{g_k}\right)\right].
\end{equation}
In estimating $\kappa(j)$, we consider the three cases $j=0,1\leq j\leq 2\nu-1$ and $2\nu\leq j\leq N-1$ separately.

\subsubsection{$j=0$}
Invoking Lemma \ref{lem:bdg}, we get a bound for $\kappa(0)$:
\[
D(N,\nu)\kappa(j)\leq \frac{3}{2}\phi(6)I_1+\frac{3}{2}\phi(4)I_2+3\phi(4)I_3=:J,
\]
where
\[
I_1=g_{2\nu-1}+g_{2\nu},\quad I_2=\sum_{k=0}^{2\nu-2}g_k+g_{k+1},\quad I_3=\sum_{k=2\nu}^{N-1}g_k+g_{k+1}.
\]
\begin{proposition}
We have for $j=0$
\begin{eqnarray*}
I_1&=&2(B_{2\nu-1}+B_{2\nu})+A_{N-2\nu+1}, \\
I_2&=& 2A_{2\nu-1}-2+A_N-A_{N-2\nu+1}+A_{N-2\nu}(A_{2\nu}-2), \\
I_3&=& 2A_N-A_{2\nu}+A_{N-2\nu}-A_{2\nu}A_{N-2\nu}-1.
\end{eqnarray*}
\end{proposition}
\begin{proof}
Insert the formulae from Proposition \ref{prop:formelgj0}, use the recurrences (\ref{eq:recinvA}),(\ref{eq:recinvB}) for $A_k$ and $B_k$ and Lemmas \ref{lem:sumid} and \ref{lem:trigid}.
\end{proof} 
With this proposition and the identity $A_N=A_{N-2\nu}A_{2\nu}+3B_{N-2\nu}B_{2\nu}$ (Lemma \ref{lem:trigid}) we see that 
\begin{eqnarray*}
J&=&\frac{3}{2}\phi(6)\left[2(B_{2\nu-1}+B_{2\nu})+A_{N-2\nu+1}\right] \\
&&+\frac{3}{2}\phi(4)\left[4A_N-4+3B_{2\nu}B_{N-2\nu}-A_{N-2\nu+1}+2(A_{2\nu-1}-A_{2\nu})\right].
\end{eqnarray*}
Now recall that $D(N,\nu)=2A_N+\frac{3}{2}B_{2\nu}B_{N-2\nu}-2$. If we then use the recurrences (\ref{eq:recinvA}), (\ref{eq:recinvB}) for $A_{2\nu-1}$ and $B_{2\nu-1}$ and set $s:=\frac{3}{2}(\phi(6)-\phi(4))=\frac{138}{1225}$ it follows with $\phi(4)=\frac{17}{25}$ that
\begin{equation}
J=\frac{51}{25}D(N,\nu)+s(6B_{2\nu}-2A_{2\nu}+A_{N-2\nu+1}).
\end{equation}
If we plug in the estimate for $B_{2\nu}$ from Lemma \ref{lem:asym} and remark that $2\nu\leq N-1$ and $N-2\nu+1\leq N-1$, we get
\begin{equation}
J\leq \frac{51}{25}D(N,\nu)+s(2\sqrt{3}-1)A_{N-1}.
\end{equation}
Using again Lemma \ref{lem:asym} on $A_{N-1}$, we obtain
\[
J\leq \frac{51}{25}D(N,\nu)+ (A_N+\sqrt{3})\frac{s}{\lambda}(2 \sqrt{3}-1).
\]
Finally, the definition of $D(N,\nu)$ and the fact that the function $\nu\mapsto B_{2\nu}B_{N-2\nu}$ is concave for $1\leq\nu\leq (N-1)/2$ and therefore attains its minimum at the border for $2\nu=N-1$ yield
\[
A_N+\sqrt{3}\leq \frac{D(N,\nu)}{2}=A_N+\frac{3}{4}B_{2\nu}B_{N-2\nu}-1\quad\text{for }N\geq 3.
\]
Thus, $\kappa(0)$ admits the bound
\[
\kappa(0)\leq \frac{51}{25}+\frac{s}{2\lambda}(2\sqrt{3}-1)\approx 2.07719\quad\text{for }N\geq 3.
\]
For $N<3$, this estimate follows from the numerical results of Table \ref{tab:numvals} on page \pageref{tab:numvals}.
\subsubsection{$1\leq j\leq 2\nu-1$}
As for $j=0$, Lemma \ref{lem:bdg} yields a bound for $\kappa(j)$ in the case $1\leq j\leq 2\nu-1$:
\[
D(N,\nu)\kappa(j)\leq \frac{3}{2}\phi(6)I_1+\frac{3}{2}\phi(4)I_2+3\phi(4)I_3=:J,
\]
where now
\[
I_1=g_0+g_1+g_{2\nu-1}+g_{2\nu},\quad I_2=\sum_{k=1}^{2\nu-2}g_k+g_{k+1},\quad I_3=\sum_{k=2\nu}^{N-1}g_k+g_{k+1}.
\]
\begin{proposition}\label{prop:formI}
We have for $1\leq j\leq 2\nu-1$
\begin{eqnarray*}
I_1&=& 2(B_j+B_{j-1}+B_{2\nu-j}+B_{2\nu-j-1})+3B_{N-2\nu+1}(B_j+B_{2\nu-j})+A_{N-j+1}+A_{N-2\nu+j+1},\\
I_2&=& 2D(N,\nu)-3(B_{2\nu-j}+B_j)(B_{N-2\nu+1}+2B_{N-2\nu})\\
&& +2(A_{j-1}-A_{N-j}+A_{2\nu-j-1}-A_{N-2\nu+j})-A_{N-j+1}-A_{N-2\nu+1+j},\\
I_3&=& A_{N-j}+A_{N-2\nu+j}-A_{2\nu-j}-A_j+3B_{N-2\nu}(B_j+B_{2\nu-j}).
\end{eqnarray*}
\end{proposition}
\begin{proof}
As in the case $j=0$, it suffices to insert the formulae from Proposition \ref{prop:formelgj0}, to use Lemmas \ref{lem:sumid} and \ref{lem:trigid} and employ the recurrences (\ref{eq:recA}),(\ref{eq:recB}),(\ref{eq:recinvA}) and (\ref{eq:recinvB}) for $A_k$ and $B_k$.
\end{proof}
Now recall that we defined $s=\frac{3}{2}(\phi(6)-\phi(4))=\frac{138}{1225}$ and  $\phi(4)=\frac{17}{25}$; thus inserting Proposition \ref{prop:formI} into the definition of $J$ and using the recursions (\ref{eq:recinvA}) and (\ref{eq:recinvB}) for $A_{j-1},B_{j-1},A_{2\nu-j-1},$ $B_{2\nu-j-1}$ yield
\begin{eqnarray*}
J&=&\frac{51}{25}D(N,\nu)+2s(3B_j-A_j+3B_{2\nu-j}-A_{2\nu-j}) \\
&&+s(A_{N-j+1}+A_{N-2\nu+j+1}+3B_{N-2\nu+1}(B_{2\nu-j}+B_j))=:J_1+J_2+J_3.
\end{eqnarray*}
From Lemma \ref{lem:asym} we deduce
\[
J_2\leq 2s(3-\sqrt{3})(B_j+B_{2\nu-j}).
\]
Since the functions $x\mapsto A_x+A_{K-x}$ and $x\mapsto B_x+B_{K-x}$ are convex for $K>0$ and $0\leq x\leq K$, we see that the maximum is attained at the border, so we get
\[
J\leq \frac{51}{25}D(N,\nu)+2s(3-\sqrt{3})(1+B_{2\nu-1})+s(A_N+A_{N-2\nu+2}+3B_{N-2\nu+1}(1+B_{2\nu-1})).
\]
We now require $\nu\geq 2$. Since we are in the case $1\leq j\leq 2\nu-1$, we see that the only case missing is $\nu=1,j=1$ which was treated above in Section \ref{sect:nu1j1}. If we now use the estimates 
\setlist{noitemsep}
\setlength{\parindent}{0.5cm}
\begin{enumerate}[labelindent=\parindent,leftmargin=*,label=\roman*.,widest=III,align=left]
\item $B_{2\nu-1}\leq \lambda^{-1}B_{2\nu}\leq \lambda^{-1}B_{2\nu}B_{N-2\nu}\quad$ (Lemma \ref{lem:asym}),
\item $A_{N-2\nu+2}\leq A_{N-2}\leq \lambda^{-2}A_N+\frac{\sqrt{3}}{\lambda}(1+\lambda^{-1})\quad$ (Lemma \ref{lem:asym}),
\item $3B_{N-2\nu+1}B_{2\nu-1}\leq \frac{A_N}{2}$\quad (Lemmas \ref{lem:asym} and \ref{lem:trigid}),
\item $3B_{N-2\nu+1}\leq 3B_{N-3}\leq 3\lambda^{-3}B_N\leq \sqrt{3}\lambda^{-3}A_N\quad$ (Lemma \ref{lem:asym}),
\end{enumerate}
we get
\begin{eqnarray*}
J-\frac{51}{25}D(N,\nu)&\leq& s(a_1+a_2A_N+a_3B_{2\nu}B_{N-2\nu})\\
&=&s\left(a_1+a_2+\frac{a_2}{2}D(N,\nu)-(\frac{3a_2}{4}-a_3)B_{2\nu}B_{N-2\nu}\right)
\end{eqnarray*}
with $a_1=2(3-\sqrt{3})+\frac{\sqrt{3}}{\lambda}(1+\lambda^{-1}),$ $a_2=\frac{3}{2}+\lambda^{-2}+\sqrt{3}\lambda^{-3}$ and $a_3=\frac{2}{\lambda}(3-\sqrt{3})$.
Since the function $\nu\mapsto B_{2\nu}B_{N-2\nu}$ is concave and therefore attains its minimum for $2\nu=N-1$ we conclude with the fact that $\frac{3a_2}{4}-a_3\geq 0$ and the exact value of this constant that
\[
J\leq D(N,\nu)\left[\frac{51}{25}+\frac{sa_2}{2}\right]\quad\text{for }N\geq 4.
\]
Thus we obtain finally
\begin{equation}
\kappa(j)\leq \frac{51}{25}+\frac{sa_2}{2}\leq 2.130411\quad\text{for }N\geq 4. \label{eq:bdnotmain}
\end{equation}
Once again, Table \ref{tab:numvals} on page \pageref{tab:numvals} yields that we have the same bound for $\kappa$ for $N<4$.
\subsubsection{$2\nu\leq j\leq N-1$}
We invoke again Lemma \ref{lem:bdg} to get
\[
D(N,\nu)\kappa(j)\leq \frac{3}{2}\phi(6)I_1+\frac{3}{2}\phi(4)I_2+3\phi(4)I_3=:J,
\]
where
\[
I_1=g_0+g_1+g_{2\nu-1}+g_{2\nu},\quad I_2=\sum_{k=1}^{2\nu-2}g_k+g_{k+1},\quad I_3=\sum_{k=2\nu}^{N-1}g_k+g_{k+1}.
\]
\begin{proposition}\label{prop:Ij2}
We have for $2\nu\leq j\leq N-1$
\begin{eqnarray*}
I_1&=& (1+B_{2\nu}+B_{2\nu-1})(A_{j-2\nu}+A_{N-j})+B_j+B_{j-1}+B_{N-j}+B_{N-j+1}\\
&&+B_{j-2\nu}+B_{j-2\nu+1}+B_{N-j+2\nu}+B_{N-j+2\nu-1},\\
I_2&=& A_{j-1}-A_{j-2\nu+1}+(A_{j-2\nu}+A_{N-j})(A_{2\nu-1}-2)+A_{N-j+2\nu-1}-A_{N-j+1},\\
I_3&=& D(N,\nu)+(1-A_{2\nu})(A_{j-2\nu}+A_{N-j})-\frac{3}{2}B_{2\nu}(B_{N-j}+B_{j-2\nu}).
\end{eqnarray*}
\end{proposition}
\begin{proof}
Insert the formulae for $g$ from Proposition \ref{prop:formelgj0} and use Lemmas \ref{lem:sumid}, \ref{lem:trigid} and the recurrences \rec for $A_k$ and $B_k$.
\end{proof}
If we apply the recurrences \rec for $A_k$ and $B_k$, Lemma \ref{lem:trigid} and Proposition \ref{prop:Ij2} to $J$, we see that it simplifies to (recall that $s=\frac{3}{2}(\phi(6)-\phi(4))=\frac{138}{1225}$ and $\phi(4)=\frac{17}{25}$)
\begin{eqnarray*}
J&=&\frac{51}{25}D(N,\nu)+s\left[3B_j-A_j+(A_{j-2\nu}+A_{N-j})(3B_{2\nu}-A_{2\nu}) \right.\\
&&\left.+3B_{N-j+2\nu}-A_{N-j+2\nu}+A_{N-j+1}+A_{j-2\nu+1}\right].
\end{eqnarray*}
Remember that $2\nu\leq j\leq N-1$. Since the functions $j\mapsto A_{N-j+1}+A_{j-2\nu+1}$, $j\mapsto 3B_j-A_j+3B_{N-j+2\nu}-A_{N-j+2\nu}$, $j\mapsto A_{j-2\nu}+A_{N-j}$ are convex, they attain their maximum at the border, in our case for $j=2\nu$, so it holds that
\[
J\leq \frac{51}{25}D(N,\nu)+s\left[6B_{2\nu}-2A_{2\nu}+3B_N-A_N+A_{N-2\nu}(3B_{2\nu}-A_{2\nu})+2+A_{N-2\nu+1}\right].
\]
For $2\nu=N-1$, we see with an estimate utilizing Lemma \ref{lem:asym} and the recurrences for $A_k$ and $B_k$ that $\kappa(j)\leq\frac{J}{D(N,\nu)}\leq \frac{51}{25}+\frac{3}{4}s\approx 2.1245$ for $N\geq 4$. If $2\nu\leq N-2$, we use the estimates 
\setlist{noitemsep}
\setlength{\parindent}{0.5cm}
\begin{enumerate}[labelindent=\parindent,leftmargin=*,label=\roman*.,widest=III,align=left]
\item $\sqrt{3}B_{2\nu}\leq A_{2\nu}\quad$ (Lemma \ref{lem:asym}),
\item $3B_N\leq \sqrt{3}A_N\quad$ (Lemma \ref{lem:asym}),
\item $A_{N-2\nu+1}\leq A_{N-1}$,
\item $A_{N-2\nu}\leq \sqrt{3}B_{N-2\nu}+1\quad$ (Lemma \ref{lem:asym}),
\item $3B_{N-2\nu}B_{2\nu}\leq A_N/2\quad$ (Lemmas \ref{lem:asym} and \ref{lem:trigid}),
\item $A_{N-1}\leq \lambda^{-1}(A_N+\sqrt{3})\quad$ (Lemma \ref{lem:asym}),
\item $B_{2\nu}\leq \frac{B_{2\nu} B_{N-2\nu}}{4}\quad (2\nu\leq N-2)$ 
\end{enumerate}
and obtain further
\begin{eqnarray*}
J-\frac{51}{25}D(N,\nu)&\leq& s\left[a_1+a_2 A_N+a_3 B_{2\nu}B_{N-2\nu}\right]\\
&=&s\left(a_1+a_2+\frac{a_2}{2}D(N,\nu)-(\frac{3a_2}{4}-a_3)B_{2\nu}B_{N-2\nu}\right)
\end{eqnarray*}
with $a_1=2+\sqrt{3}\lambda^{-1}$, $a_2=\frac{3}{2}(\sqrt{3}-1)+\lambda^{-1}$, $a_3=\frac{3}{4}(3-\sqrt{3})$. Since $\frac{3}{4}a_2-a_3>0$, we conclude that
\[
\kappa(j)=\frac{J}{D(N,\nu)}\leq \frac{51}{25}+\frac{sa_2}{2}\approx 2.117\quad\text{for }N\geq 5.
\]
For $N<5$, see Table \ref{tab:numvals} on page \pageref{tab:numvals}.
\paragraph{Summary} What we have shown up to now is that in particular for $N\geq 5$ even, for all $1\leq \nu\leq \frac{N-1}{2}$ and all $0\leq j\leq N-1$ (\emph{except} the case $\nu=1,j=1$)
\begin{equation}\label{eq:kappaend}
\kappa(j)\leq 2.130411, \quad\text{(see (\ref{eq:bdnotmain}))}.
\end{equation}
\subsubsection{$\kappa(j)$ for $N$ odd}\label{sec:kappaodd}
Now let $N$ be odd. We recall the formula (\ref{eq:kappa}) for $\kappa(j)$
\[
\kappa(j)=D(N,\nu)^{-1}\left[\frac{3}{2}\sum_{k=0}^{2\nu-1}(|g_k|+|g_{k+1}|)\cdot \xi_{j,k}+3\sum_{k=2\nu}^{N-1}(|g_k|+|g_{k+1}|)\cdot \xi_{j,k}\right],
\]
where 
\[
\xi_{j,k}=\begin{cases}
1, & \text{if }\sgn a_{j,k}=\sgn a_{j,k+1},\\
\phi(|g_{k+1}|/|g_k|),& \text{otherwise}
\end{cases}.
\]
If we write formula (\ref{eq:kappa}) in the form $\kappa(j)=\sum_{k=0}^{N-1}s_k$, every summand $s_k$ admits the (trivial) bound
\[
s_k\leq \frac{3(|g_k|+|g_{k+1}|)}{D(N,\nu)},
\]
since $\phi(t)\leq 1$ for all $t\geq 0$. We now call $D^e(N,\nu)$ and $g^e_k$ the formulae for $D(N,\nu)$ and $g_k$ respectively, but for $N$ even. That is, write $1$ instead of $(-1)^N$ in formula (\ref{eq:D}) and the expressions for $g_k$ in Proposition \ref{prop:formelgj0}, no matter if $N$ is even or odd.
Then we get further
\begin{equation}\label{eq:s1}
s_k\leq \frac{3(g^e_k+g^e_{k+1})}{D^e(N,\nu)}.
\end{equation}
Easy estimates for $g^e_k$ and $D^e(N,\nu)$ supply us now with
\begin{equation}\label{eq:s2}
\frac{3(g^e_k+g^e_{k+1})}{D^e(N,\nu)}\leq 10^{-3},
\end{equation}
provided $\frac{N-3}{2}\leq |k-j|\leq \frac{N+3}{2}$ and $N\geq 19$. So, let $N\geq 19$. Define the index set $\Lambda=\left\{\frac{N-3}{2},\frac{N-1}{2},\frac{N+1}{2},\frac{N+3}{2}\right\}$. Then
\[
\kappa(j)=\sum_{k=0}^{N-1}s_k = \sum_{k\notin\Lambda}s_k+\sum_{k\in\Lambda}s_k.
\]
We obtain further that $\sum_{k\notin\Lambda}s_k$ equals
\begin{eqnarray*}
D(N,\nu)^{-1}\left[\frac{3}{2}\sum_{\stackrel{k=0}{k\notin\Lambda}}^{2\nu-1}(|g_k|+|g_{k+1}|)\phi(|g_{k+1}|/|g_k|)+3\sum_{\stackrel{k=2\nu}{k\notin\Lambda}}^{N-1}(|g_k|+|g_{k+1}|)\phi(|g_{k+1}|/|g_k|)\right]
\end{eqnarray*}
and by the above considerations this is less or equal
\begin{eqnarray*} D^e(N,\nu)^{-1}\left[\frac{3}{2}\sum_{\stackrel{k=0}{k\notin\Lambda}}^{2\nu-1}(g_k^e+g_{k+1}^e)\phi(|g_{k+1}|/|g_k|)+3\sum_{\stackrel{k=2\nu}{k\notin\Lambda}}^{N-1}(g_k^e+g_{k+1}^e)\phi(|g_{k+1}|/|g_k|)\right]
\end{eqnarray*}
We apply Lemma \ref{lem:schrunger} and see that the terms $\phi(|g_{k+1}|/|g_k|)$ admit the same bounds as for the case $N$ even. Thus, if we first apply the estimate and then omit the restriction $k\notin\Lambda$ for the summation scope, we arrive at estimating the same sum as for the case $N$ even. Since for the case $N$ even we got the bound (\ref{eq:kappaend}) (except for $\nu=j=1$), we obtain finally
\[
\sum_{k\notin\Lambda}s_k\leq 2.130411.
\]
The remaining sum $\sum_{k\in\Lambda}s_k$ is now estimated using (\ref{eq:s1}) and (\ref{eq:s2}) and we get 
\[
\sum_{k\in\Lambda}s_k\leq 4\cdot 10^{-3},
\]
so, if we summarize, we get
\[
\kappa(j)\leq  2.134411
\]
for all $N\geq 19,\nu,j$ (no matter if $N$ is odd or even) \emph{except} the case $\nu=j=1$.
\paragraph{Summary}
Thus if we combine the present section (Section \ref{sec:estkappa}) with Sections \ref{sec:eqspaced} and \ref{sect:nu1j1}, we have now shown that for all $N\geq 19$, $0\leq \nu\leq n$ and $0\leq j\leq N-1$, we have the bound
\[
\kappa(j)<\gamma.
\]
The numerical results of Table \ref{tab:numvals} on page \pageref{tab:numvals} yield this estimate for $N\leq 20$, so we get the first assertion of our main theorem (i.e. that $\left\|P_{n,\nu}:L^{\infty}(\mathbb{T})\rightarrow L^{\infty}(\mathbb{T})\right\|<\gamma$ for all $n\in\mathbb{N},0\leq \nu\leq n$). The asymptotic value $\gamma$ for $\left\|P_{n,1}:L^{\infty}(\mathbb{T})\rightarrow L^{\infty}(\mathbb{T})\right\|$ (as $n\rightarrow\infty$) was already identified in Section \ref{sect:nu1j1}. So, the proof Theorem \ref{th:main} is complete.

\begin{landscape}
\begin{table}
\begin{small}
\begin{tabular}{|c||c|c|c|c|c|c|c|c|c|c|}
\hline

$\nu\rightarrow$ & \multirow{2}{*}{0} & \multirow{2}{*}{1} & \multirow{2}{*}{2} & \multirow{2}{*}{3} & \multirow{2}{*}{4} & \multirow{2}{*}{5} & \multirow{2}{*}{6} & \multirow{2}{*}{7} &\multirow{2}{*}{8}&\multirow{2}{*}{9}\\ 
$N \downarrow$  &  &  &  &  &  &  &  & && \\
\hline\hline

2 & 1.66666667 &&&&&&&&& \\
\hline
3 & 1.77777778 & 1.84444444 &&&&&&&&\\
\hline
4 & 1.88888889 & 2.00000000 &&&&&&&&\\
\hline
5 & 1.94696970 & 2.06951872 & 1.99530864 &&&&&&&\\
\hline
6 & 1.96835017 & 2.09951691 & 2.03615841&&&&&&&\\
\hline
7 & 1.98631436 & 2.12227384 & 2.05943912 & 2.03242817&&&&&&\\
\hline
8 & 1.99137719 & 2.12904795 & 2.06731688 & 2.05587710&&&&&&\\
\hline
9 & 1.99637151 & 2.13550178 & 2.07350359 & 2.06635304 & 2.04175181&&&&&\\
\hline
10 & 1.99767916  & 2.13721808 & 2.07535577 & 2.06916162 & 2.06184314&&&&&\\
\hline
11 & 1.99903054 & 2.13897416 & 2.07709926 & 2.07093598 & 2.06877403 & 2.04423294&&&&\\
\hline
12 & 1.99937738 & 2.13942680 & 2.07756643 & 2.07147842 & 2.07063242 & 2.06343762&&&&\\
\hline
13 & 1.99974043 & 2.13989929 & 2.07804184 & 2.07192928 & 2.07140616 & 2.06942343 & 2.04489705&&&\\
\hline
14 & 1.99983312 & 2.14002005 & 2.07816469 & 2.07206051 & 2.07160716 & 2.07106530 & 2.06386325&&&\\
\hline
15 & 1.99993046 & 2.14014679 & 2.07829271 & 2.07218617 & 2.07173511 & 2.07157865 & 2.06959951 & 2.04507495&&\\
\hline
16 & 1.99995528 & 2.14017911 & 2.07832546 & 2.07221958 & 2.07177393 & 2.07171306 & 2.07118008 & 2.06397716 && \\
\hline
17 & 1.99998137 & 2.14021308 & 2.07835981 & 2.07225375 & 2.07180634 & 2.07176873 & 2.07162518 & 2.06964688 & 2.04512262 &\\
\hline
18 & 1.99998802 & 2.14022174 & 2.07836857 & 2.07226256 & 2.07181575 & 2.07178315 & 2.07174416 & 2.07121066 & 2.06400767 &\\
\hline
19 & 1.99999501 & 2.14023084 & 2.07837778 & 2.07227176 & 2.07182477 & 2.07179234 & 2.07178110 & 2.07163781 & 2.06965959 & 2.04513539 \\
\hline
20 & 1.99999679 & 2.14023316 & 2.07838012 & 2.07227411 & 2.07182717 & 2.07179513 & 2.07179076 & 2.07175241 & 2.07121883 & 2.06401584 \\ \hline

\end{tabular}
\end{small}
\caption{Values of $\left\|P_{n,\nu}\right\|_\infty$ for different values of $\nu,N=n+\nu$ obtained with \textsc{Mathematica}. The numbers are rounded to the last digit. We have with the same precision $\gamma \approx 2.14023734$.}
\label{tab:numvals}
\end{table}

\end{landscape}

\appendix
\section{Appendix}
\begin{proof}[Proof of Lemma \ref{lem:bdg}]
In order to prove (\ref{eq:schrg1}) and (\ref{eq:schrg2}) we recall the bounds (\ref{eq:A4}) and (\ref{eq:B4}), which are
\begin{equation}
B_{l+1}\leq 4B_l\text{ for }l\geq 1, \quad A_{l+1}\leq 4A_l\text{ for }l\geq 0.\label{eq:bdAB}
\end{equation}
We consider several cases depending on the values of $j,k,\nu$:
\setlist{noitemsep}
\setlength{\parindent}{0cm}
\begin{enumerate}[labelindent=\parindent,leftmargin=*,label=\textsc{Case }\Roman*.,widest=III,align=left]
\item $0\leq j\leq 2\nu-1$
\begin{enumerate}[labelindent=-1.5cm,leftmargin=*,label=\theenumi\alph*.,align=left,widest=I]
\item $k=0,j\neq 0$\\
If we note (\ref{eq:bdAB}) and the formula for $g_k$ from Proposition \ref{prop:formelgj0}, we get the inequality $6g_1-g_0\geq 0$ immediatly. For the reversed one we get, since we assumed $j\geq 1$
\begin{eqnarray*}
6g_0-g_1&=&(12B_j-2B_{j-1})+(6B_{N-j}-B_{N-j+1})\\
&&+B_{2\nu-j}(6A_{N-2\nu}-A_{N-2\nu+1})-A_{N-j}-3B_{2\nu-j}B_{N-2\nu}\\
&\geq & 10B_j+2B_{N-j}+2B_{2\nu-j}A_{N-2\nu}-A_{N-j}-3B_{2\nu-j}B_{N-2\nu},
\end{eqnarray*}
by (\ref{eq:bdAB}). If we now additionally observe that $2B_{N-j}\geq A_{N-j}$ (for $N-j\geq 1$, which is satisfied) and $A_{N-2\nu}\geq \sqrt{3}B_{N-2\nu}$, we see that this is $\geq 0$.

\item $1\leq k\leq j-1$ \\
Again, with (\ref{eq:bdAB}) and the assumption $k\leq j-1$ we get the first inequality $4g_k-g_{k+1}\geq 0$ immediatly. The second inequality is only critical for $k=j-1$ and in this case we get (with (\ref{eq:bdAB}))
\[
4g_{k+1}-g_k=-2+4B_N-B_{N-1}+\text{(positive term)}\geq 3B_N-2\geq 0\quad\text{for }N\geq 1.
\]
\item $j\leq k\leq 2\nu-2$\label{it:caseic}\\
For the first inequality $4g_{k+1}-g_k\geq 0$, it suffices to argue with (\ref{eq:bdAB}), so it does for the second one $4g_{k}-g_{k+1}\geq 0$ in the case $k\neq j$. For $k=j$ it holds that
\[
4g_k-g_{k+1}=-2+4B_{N}-B_{N-1}+\text{(positive term)}\geq 0\quad\text{for }N\geq 1.
\]
\item $k=2\nu-1$ \\
An analogous distinction between the cases $k=j$ and $k>j$ as in  \ref{it:caseic} supplies us with the estimate $6g_k-g_{k+1}\geq 0$. On the other hand (recall that $k=2\nu-1$,$j\leq 2\nu-1$)
\begin{eqnarray*}
6g_{k+1}-g_k&=&(12B_{2\nu-j}-2B_{2\nu-1-j})+(6B_{N-2\nu+j}-B_{N-2\nu+1+j})\\
&& -A_{N-2\nu+j}+(6B_jA_{N-2\nu}-B_jA_{N-2\nu+1})-3B_jB_{N-2\nu} \\
&\geq & 0+2B_{N-2\nu+j}-A_{N-2\nu+j}+2B_jA_{N-2\nu}-3B_jB_{N-2\nu},
\end{eqnarray*}
by (\ref{eq:bdAB}). The inequalities $2B_{N-2\nu+j}\geq A_{N-2\nu+j}$ (observe $N-2\nu+j\geq 1$) and $A_{N-2\nu}\geq \sqrt{3}B_{N-2\nu}$ then yield $6g_{k+1}-g_k\geq 0$.
\item $2\nu\leq k\leq N-1$\\
Since $k>j$ in the current case, an application of (\ref{eq:bdAB}) suffices for $4g_k-g_{k+1}\geq 0$. The same reasoning provides us with $4g_{k+1}-g_k\geq 0$ in the case $k\neq N-1\vee j\neq 0$ and for $k=N-1,j=0$ we have
\begin{eqnarray*}
4g_{k+1}-g_k&=&(4B_N-B_{N-1})+B_{2\nu}(4A_{N-2\nu}-A_{N-1-2\nu})-1\\
&\geq & 3B_N-1\geq 0 \quad\text{for }N\geq 1.
\end{eqnarray*}
\end{enumerate}
\item $2\nu\leq j\leq N-1$

\begin{enumerate}[labelindent=-1.5cm,leftmargin=*,label=\theenumi\alph*.,align=left,widest=I]
\item $k=0$\\
Again, the estimate $6g_{k+1}-g_k\geq 0$ is a trivial consequence of (\ref{eq:bdAB}). Furthermore, by (\ref{eq:bdAB}),
\[
6g_k-g_{k+1}=6B_{N-j}-B_{N-j+1}-A_{N-j}+\text{(positive terms)}\geq 2B_{N-j}-A_{N-j}\geq 0.
\]
\item $1\leq k\leq 2\nu-2$ \\
Here, both inequalities $4g_{k}-g_{k+1}\geq 0$ and $4g_{k+1}-g_k\geq 0$ are a consequence of (\ref{eq:bdAB}).
\item $k=2\nu-1$ \\
The bound $6g_{k}-g_{k+1}\geq 0$ follows from (\ref{eq:bdAB}). For the converse we get
\[
6g_{k+1}-g_k=(6B_{j-2\nu}-B_{j-2\nu+1})-A_{j-2\nu}+(6B_{N-j+2\nu}-B_{N-j+2\nu-1})
+\text{(positive term)}.
\]
If $j>2\nu$, we have $6B_{j-2\nu}-B_{j-2\nu+1}\geq 2B_{j-2\nu}$, which is greater than $A_{j-2\nu}$; if $j=2\nu$, $6g_{k+1}-g_k\geq -2+5B_N\geq 0.$
\item $2\nu\leq k\leq j-1$\label{it:IId}\\
For $k>2\nu$, $4g_k-g_{k+1}\geq 0$ is a consequence of (\ref{eq:bdAB}). If $k=2\nu$, we have
\[
4g_k-g_{k+1}= 2B_{N-j+2\nu}-\frac{3}{2}B_{2\nu}B_{N-j}+\text{(positive term)}.
\]
Since $2B_{N-j+2\nu}\geq A_{N-j+2\nu}$ and $3B_{2\nu}B_{N-j}=A_{N-j+2\nu}-A_{2\nu}A_{N-j}\leq A_{N-j+2\nu}$, we get $4g_k-g_{k+1}\geq 0$. The converse estimate $4g_{k+1}-g_k\geq 0$ follows once more from (\ref{eq:bdAB}) provided $k< j-1$. If on the other hand we have $k=j-1$, we see that
\[
4g_{k+1}-g_k=-1+A_{N-j}(4B_{k+1}-B_k)+\text{(positive term)}\geq 0,
\]
since $k=j-1\geq 2\nu\geq 2$.
\item $j\leq k\leq N-1$ \\
The estimate $4g_k-g_{k+1}\geq 0$ follows from (\ref{eq:bdAB}) if $k>j$, as does $4g_{k+1}-g_k\geq 0$ for $k<N-1$. For the critical values $k=j$ resp. $k=N-1$, similar calculations as in \ref{it:IId} conclude the statement of the lemma.\qedhere
\end{enumerate}
\end{enumerate}

\end{proof}

\begin{proof}[Proof of Lemma \ref{lem:schrunger}]
If $N$ is odd, the proof consists of similar estimates as the proof of Lemma \ref{lem:bdg} but with twice as many case distinctions, since one has to consider the cases $|k-j|\leq\frac{N-5}{2}$ and $|k-j|\geq \frac{N+5}{2}$ separately. We pick out one special case and omit all the others since they involve very similar arguments to the presented case or even to the proof of Lemma \ref{lem:bdg}. 
We will treat values of $\nu,k,j$ where $2\nu\leq k\leq j-1$ and view the two cases mentioned above:
\setlist{noitemsep}
\setlength{\parindent}{0cm}
\begin{enumerate}[labelindent=\parindent,leftmargin=*,label=\textsc{Case }\Roman*.,widest=III,align=left]
\item $|j-k|\leq \frac{N-5}{2}$ \\
We obtain from Proposition \ref{prop:formelgj0} and Remark \ref{eq:remgpos} that
\begin{eqnarray*}
|g_k|&=&-B_{j-k}+A_{k-2\nu}B_{N-j+2\nu}+A_{N-j}B_k+\frac{3}{2}B_{k-2\nu}B_{2\nu}B_{N-j}, \\
|g_{k+1}|&=&-B_{j-k-1}+A_{k+1-2\nu}B_{N-j+2\nu}+A_{N-j}B_{k+1}+\frac{3}{2}B_{k+1-2\nu}B_{2\nu}B_{N-j}.
\end{eqnarray*}
The inequality $4|g_k|-|g_{k+1}|\geq 0$ for $k=2\nu$ is a simple consequence of Lemmas \ref{lem:asym} and \ref{lem:trigid}. Utilizing Lemma \ref{lem:asym}, we get for $k\geq 2\nu+1$ that
\begin{equation}
4|g_k|-|g_{k+1}|\geq -4B_{j-k}+(4-\lambda)A_{k-2\nu}B_{N-j+2\nu}. \label{eq:estungerade}
\end{equation}
Since $N-j+2\nu\geq 3$, $A_3=26$ and $2\nu\leq k$ we see with Lemma \ref{lem:trigid} that
\[
A_{k-2\nu}\leq \frac{A_{k-2\nu} A_{N-j+2\nu}}{A_3}= \frac{A_{k-2\nu}A_{N-j+2\nu}}{26}\leq \frac{A_{N-j+k}}{26}.
\]
This estimate, the definition of the recurrences $A_k$ and $B_k$ and Lemmas \ref{lem:asym} and \ref{lem:trigid} yield
\begin{eqnarray*}
A_{k-2\nu}B_{N-j+2\nu}&\geq& \frac{1}{\sqrt{3}}(A_{k-2\nu}A_{N-j+2\nu}-A_{k-2\nu})\geq \frac{1}{2\sqrt{3}}(A_{N-j+k}-2A_{k-2\nu}) \\
&\geq& \frac{2\sqrt{3}}{13}A_{N-j+k}.
\end{eqnarray*}
Thus, this estimate and (\ref{eq:estungerade}) imply 
\begin{eqnarray*}
4|g_k|-|g_{k+1}|&\geq& (4-\lambda)\frac{2\sqrt{3}}{13}A_{N-j+k}-4B_{j-k}\geq (4-\lambda)\frac{6}{13}B_{N-j+k}-4B_{j-k} \\
&\geq& (\lambda^5(4-\lambda)\frac{6}{13}-4)B_{(N-5)/2}\geq 0,
\end{eqnarray*}
if we use Lemma \ref{lem:asym} in conjunction with our hypothesis $|j-k|\leq \frac{N-5}{2}$. The estimate $4|g_{k+1}|-|g_k|\geq 0$ follows analogously.
\item $|j-k|\geq \frac{N+5}{2}$ \\
We obtain from Proposition \ref{prop:formelgj0} and Remark \ref{eq:remgpos} that
\begin{eqnarray*}
|g_k|&=&B_{j-k}-A_{k-2\nu}B_{N-j+2\nu}-A_{N-j}B_k-\frac{3}{2} B_{k-2\nu}B_{2\nu}B_{N-j}, \\
|g_{k+1}|&=&B_{j-k-1}-A_{k+1-2\nu}B_{N-j+2\nu}-A_{N-j}B_{k+1}-\frac{3}{2}B_{k+1-2\nu}B_{2\nu}B_{N-j}.
\end{eqnarray*}
If we employ Lemma \ref{lem:asym} three times, we obtain
\[
4|g_{k}|-|g_{k+1}|\geq 3B_{j-k}-(4-\lambda)[B_{N-j+2\nu}A_{k-2\nu}+B_k A_{N-j}+\frac{3}{2}B_{2\nu}B_{N-j}B_{k-2\nu}]
\]
Since by Lemma \ref{lem:trigid} every summand in the square bracket is majorized by $B_{N-j+k}$, we finally get
\[
4|g_{k}|-|g_{k+1}|\geq 3(B_{j-k}-(4-\lambda)B_{N-j+k})\geq 0,
\]
by the hypothesis $|j-k|\geq \frac{N+5}{2}$.
For the inequality $4|g_{k+1}|-|g_k|\geq 0$, we first omit some positive terms to get
\[
4|g_{k+1}|-|g_k|\geq 4B_{j-k-1}-B_{j-k}-4A_{k+1-2\nu}B_{N-j+2\nu}-4A_{N-j}B_{k+1}-6B_{k+1-2\nu}B_{2\nu}B_{N-j}.
\]
As above, Lemmas \ref{lem:trigid} and \ref{lem:asym} respectively yield
\begin{eqnarray*}
4|g_{k+1}|-|g_k|&\geq &4B_{j-k-1}-B_{j-k}-10B_{N-j+k+1} \\
&\geq& (4-\lambda)B_{j-k-1}-1-10B_{N-j+k+1}.
\end{eqnarray*}
But now we employ again Lemma \ref{lem:asym} and the fact that $|j-k|\geq \frac{N+5}{2}$ to get
\[
4|g_{k+1}|-|g_k|\geq (\lambda^3(4-\lambda)-10)B_{(N-3)/2}-1\geq 0,
\]
and so the desired inequality.\qedhere
\end{enumerate}
\end{proof}

\newpage
\bibliographystyle{plain}
\bibliography{LebesgueTorus}

\nocite{Ciesielski1963}\nocite{Ciesielski1966}\nocite{CiesielskiKamont2004}\nocite{Bechler2003}
\nocite{Ciesielski1975}\nocite{Ciesielski1975a}\nocite{Oskolkov1979}\nocite{Oswald1977}\nocite{Domsta1972}\nocite{Domsta1976}\nocite{Keryan2005}\nocite{Keryan2008}\nocite{CiesielskiNiedzwiecka}\nocite{Shadrin2001}\nocite{KashinSaakyan1989}\nocite{deBoor1968}

\small
\begin{tabular}{p{10cm}l}
  \textsc{Markus Passenbrunner} \\
  \textsc{Department of Analysis} \\ 
  \textsc{J. Kepler University}\\
  \textsc{Altenberger Strasse 69}\\
  \textsc{A-4040 Linz}\\
  \textsc{Austria}\\
  passenbr@bayou.uni-linz.ac.at
\end{tabular}

\end{document}